\documentclass[10pt,reqno]{amsart}
\usepackage{amssymb}
\usepackage{amsmath}
\usepackage{bm}
\usepackage{amsthm}
\usepackage{amsbsy}
\usepackage{psfrag}
\usepackage{graphics}
\usepackage{graphicx}
\usepackage{pgfplots}
\usepackage{relsize}
\usepackage{color}
\usepackage{xcolor}
\usepackage{pdfpages}  
\usepackage{url}
\usepackage{hyperref}
\usepackage[square,sort&compress,comma,numbers]{natbib}
\usepackage{graphicx}
\usepackage{caption}
\usepackage{subcaption}
\hypersetup{
	colorlinks=true,
	linkcolor=blue,
	filecolor=magenta,      
	urlcolor=cyan,
}
\theoremstyle{plain}
\newtheorem{theorem}{Theorem}[section]
\newtheorem{proposition}[theorem]{Proposition}

\newtheorem{example}[theorem]{Example}
\newcommand{\norm}[1]{\left\Vert #1 \right\Vert}
\newtheorem{remark}[theorem]{Remark}

\begin{document}
\allowdisplaybreaks[4]
\numberwithin{figure}{section}

%
\numberwithin{table}{section}
 \numberwithin{equation}{section}
%

\title[FEM for Dirichlet Boundary Optimal Control  Problem]
{A two-grid Adaptive Finite Element Method for the Dirichlet Boundary Control Problem Governed by Stokes Equation}

\author[T. Gudi]{ Thirupathi Gudi}
 \address{Department of Mathematics, Indian Institute of Science, Bangalore - 560012, India}
 \email{gudi@iisc.ac.in}

 \author[R. C. Sau]{Ramesh Ch. Sau}

\address{Department of Mathematics, Chinese University of Hong Kong, Hong Kong}

\email{rcsau1994@gmail.com}

\date{}

\begin{abstract}
  In this article, we derive \textit{a posteriori} error estimates for the Dirichlet boundary control problem governed by Stokes equation. An energy-based method has been deployed to solve the Dirichlet boundary control problem. We employ an inf-sup stable finite element discretization scheme by using $\mathbf{P}_1$ elements(in the fine mesh) for the velocity and control variable and $P_0$ elements(in the coarse mesh) for the pressure variable. We derive an \textit{a posteriori} error estimator for the state, adjoint state, and control error. The control error estimator generalizes the standard residual type estimator of the unconstrained Dirichlet boundary control problems, by additional terms at the contact boundary addressing the non-linearity. We prove the reliability and efficiency of the estimator. Theoretical results are illustrated by some numerical experiments.
\end{abstract}


\maketitle
\allowdisplaybreaks
\def\R{\mathbb{R}}
\def\I{\mathbb{I}}
\def\dx{{\rm~dx}}
\def\dy{{\rm~dy}}
\def\ds{{\rm~ds}}
\def\y{\mathbf{y}}
\def\u{\mathbf{u}}
\def\V{\mathbf{V}}
\def\w{\mathbf{w}}
\def\J{\mathbf{J}}
\def\f{\mathbf{f}}
\def\g{\mathbf{g}}
\def\div{{\rm div~}}
\def\v{\mathbf{v}}
\def\z{\mathbf{z}}
\def\b{\mathbf{b}}
\def\cV{\mathcal{V}}
\def\V{\mathbf{V}}
\def\cT{\mathcal{T}}
\def\sjump#1{[\hskip -1.5pt[#1]\hskip -1.5pt]}

\section{Introduction}\label{sec:Intro}


The study of the optimal control problems governed by partial differential equations has been a major research area in applied mathematics and its allied areas. The optimal control problem consists of finding an optimal control variable that minimizes a cost functional subject to a partial differential equation satisfied by the optimal control and an optimal state. There are many results on the finite element analysis of optimal control problems, see for example \cite{Hinze:2011:Control,CasasM:2008:Neum,Hinze:2005:Control,Casas_Dhamo_2012}. In an optimal control problem, the control can act on the system through either a boundary
condition or through an interior force. In the latter case, the control is distributed
and in the former case, the control is said to be the boundary control. The choice of the boundary condition leads to several types of boundary controls e.g., Dirichlet, Neumann or Robin boundary control. We refer to \cite{Hinze:2011:Control,meyerrosch:2004:Optiaml,Hinze:2005:Control,Sudipto:2015:IP,Ramesh:2018} for the work related to distributed control and to \cite{CasasM:2008:Neum,Casas_Dhamo_2012,Sudipto:2015:IP,Ramesh:2018} for the work on Neumann boundary control problem.

The Dirichlet boundary control problems are important in application areas, but the problems can
be difficult to analyze mathematically because in some cases the control itself solves some PDE and in other cases the physical domain has a non-smooth boundary. The study of Dirichlet control problems posed on polygonal
domains can be traced back to \cite{CR:2013:Dirich}, where a control constrained problem governed
by a semilinear elliptic equation posed in a convex polygonal domain is studied.
There are various approaches proposed in the literature. One of the popular approach is to seek control from the $L^2(\Gamma)$-space. In this case, the state equation has to be understood in a ultra weak sense, since the Dirichlet boundary condition is only in $L^2(\Gamma)$. This ultra weak formulation is easy to implement
and usually results in optimal controls with low regularity. Especially, when the problem is posed on a convex
polygonal domain, the control $y$ vanishes on the corners and is thus continuous. This is because, it is determined
by the normal derivative of the adjoint state, whereas in a nonconvex polygonal domain the control may have singularity around the corner point, for more details one can see \cite{apel:2015}. An other approach, as in \cite{CMR:2009:Dirich}, is the Robin boundary penalization which transforms the
Dirichlet control problem into a Robin boundary control problem.

One other popular approach is to find controls from the energy space, i.e., $H^{1/2}(\Gamma),$ we refer \cite{Steinbach:2014:Dirichlet} for this approach. In \cite{Steinbach:2014:Dirichlet},
 the Steklov-Poincar\'e operator was used to define the cost functional with the help of a harmonic extension of the given boundary data. The Steklov-Poincar\'e operator transforms the Dirichlet data into a Neumann data by using harmonic extension of the Dirichlet data; 
  but this type of abstract operator may cause some difficulties in numerical implementation. A boundary element method for this problem is proposed and analyzed in \cite{Steinbach:2014}.

Given any $y$ in the Sobolev space $H^{1/2}(\Gamma)$, we can construct a function $u_y$ in the Sobolev space $H^1(\Omega)$ that is harmonic in $\Omega$ and agrees with $y$ on the boundary $\Gamma$. This function $u_y$ is called a harmonic extension of $y$. The norm of $y$ in $H^{1/2}(\Gamma)$ can be equivalently written as the norm of the gradient of $u_y$ in $L^2(\Omega)$, i.e., $|y|_{H^{1/2}(\Gamma)}\equiv \norm{\nabla u_y}_{0,\Omega}.$ This relation suggests that we can penalize the control $y$ in the $H^1(\Omega)$ space by adding the term $\norm{\nabla y}_{0,\Omega}^2$ to the cost functional. Therefore, we consider the following optimization problem:
  
  \begin{equation*}
  \text{min} ~J(u,y)= \frac{1}{2}\norm{u-u_d}_{0,\Omega}^2 + \frac{\lambda}{2}\norm{\nabla y}_{0,\Omega}^2.
  \end{equation*}

  The paper \cite{Gudi:2017:DiriControl} proposed a novel method for the Dirichlet boundary control problem using the above cost functional. The method is based on the energy space, where the control belongs to the Sobolev space $H^1(\Omega)$. This avoids the use of the Steklov-Poincar\'e operator and makes the method computationally efficient. The paper \cite{Gudi:2017:DiriControl} only considered the unconstrained case, while the constrained case was analyzed in \cite{Gudi_Ramesh:ESIAM_COCV}. The paper \cite{karkulik} presented a similar method based on the energy space. The paper \cite{Winkler} obtained a sharp convergence rate for the energy space method.
\par

 The literature on Stokes Dirichlet control problem considers two types of control. One is the tangential control, i.e., the control acts only in the tangential direction of the boundary (see \cite{gong:HDG}). The paper \cite{gong:HDG} used a hybridized discontinuous Galerkin (HDG) method to solve a tangential Dirichlet boundary control problem with an $L^2$ penalty on the boundary control, without any constraints on the control. The other is the zero flux control, which means the control has no normal component to the boundary (i.e., $\int_{\partial\Omega}\mathbf{y}\cdot\mathbf{n}=0$) \cite{gong_mateoes_stokes}. The paper \cite{gong_mateoes_stokes} studied two different boundary control regularization terms in the cost functionals: the $L^2$ norm and the energy space seminorm. The zero flux condition is a natural consequence of the incompressibility condition and the Dirichlet boundary condition in the PDE. Therefore, the authors either chose the tangential control as the first option or imposed the zero flux condition as a constraint in the space as the second option. Many papers on Navier-Stokes Dirichlet control problem also used either tangential control or zero flux control, for example, see \cite{Fursikov,Gunzburger}. The zero flux condition on the control affects the regularity of the control as discussed in \cite{gong_mateoes_stokes}. To address this issue, we introduced the Stokes equation with mixed boundary conditions and the control acts only on the Dirichlet boundary. Our control is more general and has both tangential and normal components. We also added constraints on the control. As a result, the optimal control satisfies a simplified Signorini problem.

 In this article, we propose, analyze, and test a new a posteriori error estimator for the control error. In order to derive and analyze the error estimator for the control variable, we adopt the framework presented
in \cite{Andreas_Veeser_Signorini:2015} for Signorini problem, because the control satisfies the Signorini problem. The discrete
problem consists of a discrete variational inequality for the approximate control variable and the estimator is designed for controlling its energy error. The estimator reduces
to the standard residual estimator for elliptic problem, if no contact occurs. The contributions by the estimator addressing the nonlinearity are related to the contact stresses, the
complementarity condition. 
We prove reliability and efficiency of the estimator and ensuring the equivalence with the error up to oscillation terms. 
A key ingredient of this approach is the so-called Galerkin functional. It is a modification of the residual with
respect to the corresponding linear problem with the help of a suitable approximation of the Lagrange multiplier and thus, may be seen as the residual of a linear auxiliary problem. The correction by the Lagrange multiplier is crucial for sharpness of the upper bounds in the actual contact regions. The theoretical results are
corroborated by a variety of numerical tests

\par
The rest of the article is organized as follows. In Section
\ref{sec:ModelProblem}, we formulate the Dirichlet boundary
control problem with pointwise control constraints. Therein, we discuss the well-posedness of the model problem and present the optimality system.
In Section \ref{Discrete Problems}, we define the discrete
control problem and present the discrete optimality system. In Section \ref{Sto:DC_error_analysis}, we derive \textit{a posteriori} error estimates with the help of some reconstruction solution. Section
\ref{sec:Numerics} is devoted to the numerical experiments.

\section{Continuous Problem}\label{sec:ModelProblem}
We proceed over the precise formulation of the optimization problem in brief in this section. We need the following definitions and notations before we can begin the analysis:
\subsection{Notation}\label{notation}
 Let $\Omega\subset\R^2,$ be a bounded polygonal domain, with
boundary $\partial\Omega$ consists of three non-overlapping open subsets $\Gamma_D, \;  \Gamma_C$ and $\Gamma_N$  with $\partial\Omega=\Gamma_C\cup\bar \Gamma_D\cup\Gamma_N$. The one-dimensional measure of $\Gamma_C$ is positive.  We denote any function and any space in bold notation can be understood in the vector form e.g.,
$\mathbf{x}:=(x_{1},x_{2}),$  $\mathbf{L}^2(\Omega):=[L^2(\Omega)]^2$ and $\mathbf{H}^1(\Omega):=[H^1(\Omega)]^2.$ The norm and inner product on those spaces are defined component wise. Here and throughout, the $\mathbf{L}^2(\Omega)$ norm is denoted by
$\|\cdot\|_{0,\Omega}$ and $\|\cdot\|_{k,\Omega} (k> 0)$ denotes the standard norm
on the Sobolev space $\mathbf{H}^k(\Omega)$, see for example \cite{Ciarlet:1978:FEM}. 
The trace of a vector valued function $\mathbf{x}\in \mathbf{H}^1(\Omega)$ is defined to be $\bm{\gamma}_{0}(\mathbf{x}):=(\gamma_{0}(x_1),\gamma_{0}(x_2)),$ where $\gamma_0:H^1(\Omega)\rightarrow L^2(\Gamma)$ is the
trace operator. Let  $\mathbf{x}$ and $\mathbf{y}$ are two functions, we say that  $\mathbf{x}\leq \mathbf{y}$ iff $x_1\leq y_1$ and $x_2\leq y_2$ almost everywhere in $\Omega.$

\subsection{Dirichlet Control Problem} We consider the following constrained Dirichlet boundary control problem(in energy form \cite[Section 2]{Gudi:2017:DiriControl}) governed by the Stokes equation
 \begin{equation}\label{min:j}
 \text{min} ~J(\mathbf{u},\mathbf{y})= \frac{1}{2}\norm{\mathbf{u}-\mathbf{u}_d}_{0,\Omega}^2 + \frac{\rho}{2}\norm{\nabla\mathbf{y}}_{0,\Omega}^2
 \end{equation}	
subject to,
\begin{equation}\label{stokes_eq_intro}
\begin{split}
-\Delta \mathbf{u}+\nabla{p}&=\mathbf{f} \quad \text{in}\;\Omega, \\
\nabla\cdot{\mathbf{u}}&=0 \quad \text{in}\; \Omega,\\
\bf{u}&=\mathbf{y} \quad \text{on}\; \Gamma_C,\\
\bf{u}&=\mathbf{0} \quad \text{on}\; \Gamma_D,\\
\frac{\partial \mathbf{u}}{\partial \mathbf{n}}-p \mathbf{n}&=\mathbf{0}\;\; \text{on} \;\; \Gamma_N, 
\end{split}	
\end{equation}

with the control $\mathbf{y}$ comes from the following constrained set
\begin{align*}
   \mathbf{Q}_{ad}:=\{\mathbf{x} \in \mathbf{Q}: \mathbf{y}_a\leq\bm{\gamma}_{0}(\mathbf{x})\leq \mathbf{y}_b\text{ a.e. on } \Gamma_C\}.
\end{align*}
The interior force $\mathbf{f}\in \mathbf{L}^2(\Omega),$ the regularization parameter $\rho>0,$ and $\mathbf{u}_d\in \mathbf{L}^2(\Omega)$ and the space $\mathbf{Q} := \{\mathbf{x} \in \mathbf{H}^1(\Omega): \bm{\gamma}_{0}(\mathbf{x})=\mathbf{0}\text{ on } \Gamma_D\cup \Gamma_N \}.$ The constant vectors $\mathbf{y}_a=(y^1_a,y^2_a)$, and $\mathbf{y}_b=(y^1_b,y^2_b)\in\R^2$ satisfying $y^1_a <y^2_a$ and $y^1_b <y^2_b$. Furthermore whenever $\Gamma_D$ is non empty, we assume for consistency that $y^1_a,y^1_b\leq 0$ and $y^2_a,y^2_b\geq0$ so that, the admissible set $\mathbf{Q}_{ad}$ is nonempty.

A proof of the existence of the unique solution of the control problem \eqref{min:j} can be found in \cite[Theorem 2.2]{Gudi_Ramesh_CaMwA:2023}. The following proposition states the first-order optimality system, a details can be found in \cite[Proposition 2.3]{Gudi_Ramesh_CaMwA:2023}.

\begin{proposition}\label{prop:Wellposed-C}
	There exists a unique solution $(\mathbf{u},p,\mathbf{y})\in \mathbf{H}_{D}^1(\Omega)\times L^2(\Omega)\times \mathbf{Q}_{ad}$ for the
	Dirichlet control problem $\eqref{min:j}$ and there exists an adjoint state $(\bm{\phi},r)\in \mathbf{V}\times L^2(\Omega)$ satisfying
	\begin{subequations}\label{conti_optimality_system}
	\begin{align}
	\mathbf{u}&=\mathbf{w}+\mathbf{y},\quad \mathbf{w}\in \mathbf{V},\label{eq:state_0a}\\
	a(\mathbf{w},\mathbf{z})+b(\mathbf{z},p) &=( \mathbf{f},\mathbf{z})-a(\mathbf{y},\mathbf{z}) \;\;\;
	{\rm for~all}\;\mathbf{z} \in \mathbf{V},\label{eq:state_a}\\ 
	b(\mathbf{u},q)&=0\;\quad {\rm for~all}\; q \in  L^2(\Omega),\label{eq:state_b}\\
	a(\mathbf{z},\bm{\phi})-b(\mathbf{z},r) &=( \mathbf{u-u_d},\mathbf{z}) \;\;\;
	{\rm for~all}\;\mathbf{z} \in \mathbf{V},\label{eq:adjstate_a}\\ 
	b(\bm{\phi},q)&=0\;\quad {\rm for~all}\; q\in  L^2(\Omega),\label{eq:adjstate_b}\\
	\rho\, a(\mathbf{y},\mathbf{x}-\mathbf{y})&\geq a(\mathbf{x}-\mathbf{y},\bm{\phi})-b(\mathbf{x}-\mathbf{y},r)-(\mathbf{u-u_d},\mathbf{x}-\mathbf{y})\quad\forall \mathbf{x}\in
	\mathbf{Q}_{ad},\label{eq:VI}
	\end{align}
\end{subequations}
where
 $a(\mathbf{w},\mathbf{z})=  \int_{\Omega} \nabla{\mathbf{w}}:\nabla{\mathbf{z}} \dx$ , $b(\mathbf{z},p)= - \int_{\Omega} p\nabla\cdot{\mathbf{z}} \dx,$ and the matrix product $A:B := \sum_{i,j=1}^{n}a_{ij}b_{ij}$ when $A=(a_{ij})_{1\leq i,j \leq n}$ and $B=(b_{ij})_{1\leq i,j \leq n}.$ 

\end{proposition}

\begin{remark}\label{signorini_L2_f}	
It is not hard to show from the equation \eqref{eq:VI}, that the optimal control $\mathbf{y}$ is the weak solution of the following simplified Signorini problem:
 	\begin{subequations}
 	\begin{align*}
 		-\rho\Delta \mathbf{y}&=\mathbf{0}\quad \text{in}\quad \Omega,\\
 		\mathbf{y}&=\mathbf{0}\quad \text{on}\quad \Gamma_D\cup \Gamma_N,\\
 		\mathbf{y}_a \leq \bm{\gamma}_0(\mathbf{y})&\leq \mathbf{y}_b\; \text{ a.e. on } \Gamma_C
 	\end{align*}
 	further, the following holds for almost every $x\in \Gamma_C$:
 	\begin{align*}
 		\text{if}\;  \mathbf{y}_a<\mathbf{y}(x) <\mathbf{y}_b \quad \text{then} \quad \big(\hat{\sigma}(\mathbf{y})\big)(x)&=\mathbf{0},\\
 		\text{if}\; \mathbf{y}_a\leq \mathbf{y}(x)<\mathbf{y}_b \quad \text{then} \quad \big(\hat{\sigma}(\mathbf{y})\big)(x)&\geq \mathbf{0},\\
 		\text{if}\;  \mathbf{y}_a<\mathbf{y}(x)\leq \mathbf{y}_b \quad \text{then} \quad \big(\hat{\sigma}(\mathbf{y})\big)(x)&\leq \mathbf{0},
 	\end{align*}
 \end{subequations}
 where $\hat{\sigma}(\mathbf{y}):=\rho \frac{\partial \mathbf{y}}{\partial\mathbf{n}} - \frac{\partial\bm{\phi}}{\partial\mathbf{n}}-r\mathbf{n}.$
\end{remark}

\section{Discrete Problem}\label{Discrete Problems}
\begin{figure}[!hh]
	\begin{center}
		\setlength{\unitlength}{0.7cm}
		\begin{picture}(8,6)
		\put(2,1){\line(1,0){8}} \put(2,1){\line(-1,1){4}}
		\put(6,5){\line(-1,0){8}} \put(10,1){\line(-1,1){4}}
		\thicklines\put(2,1){\vector(1,1){4}}\thinlines
		\put(1.39,0.39){$A$} \put(5.95,5.25){$B$} \put(10,0.29){$P_{+}$}
		\put(-2.75,4.85){$P_{-}$} \put(1.25,3.25){$T_{-}$}
		\put(6,2.5){$T_{+}$} \put(4,3){\vector(1,-1){1}}
		\put(5.10,1.90){$\mathbf{n}_{e}$}
		\thicklines\put(4,3){\vector(1,1){1}}\thinlines
		\put(5.10,3.50){$\tau_{e}$} \put(3.25,2.75){$e$}
		\qbezier(4.2,2.8)(4.4,3.0)(4.2,3.2) \put(4.1,2.95){.}
		\end{picture}
		\caption{\footnotesize{Here $T_-$ and $T_{+}$ are the two neighboring triangles that share the edge $e=\partial T_-\cap\partial T_+$ with initial node $A$
				and end node $B$ and unit normal $\mathbf{n}_e$. The orientation of $\mathbf{n}_e = \mathbf{n}_{-} = -\mathbf{n}_{+}$ equals the outer normal of $T_{-}$, and hence, points into $T_{+}$.}} \label{Fig1}
	\end{center}
\end{figure}
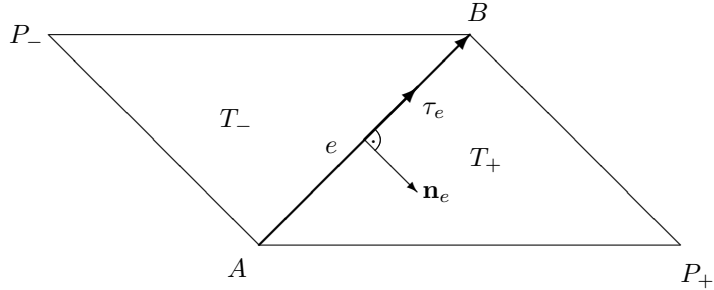

In this section, we discuss the discrete control problem before this we need to define some notations.
Let $\cT_H$ be a shape-regular triangulation of the domain $\Omega$ into triangles such that $\cup_{T\in \cT_H}T=\bar{\Omega}$ see \cite{BScott:2008:FEM,Ciarlet:1978:FEM}. Also let $\cT_h$ be a refinement of $\cT_H$ by connecting all the midpoints of $\cT_H.$ Denote the set of all interior edges of $\cT_h$ by $\mathcal{E}_{h}^{i}.$ The set of all Dirichlet, Neumann and Contact boundary edges of $\cT_h$ are denoted by $\mathcal{E}_{h}^{b,D},$ $\mathcal{E}_{h}^{b,N}$ and $\mathcal{E}_{h}^{b,C}$ respectively and define $\mathcal{E}_{h}=\mathcal{E}_{h}^{i}\cup\mathcal{E}_{h}^{b,D}\cup \mathcal{E}_{h}^{b,N}\cup\mathcal{E}_{h}^{b,C}.$ A typical triangle is denoted by $T$ and its
diameter by $h_T$. Set $h=\max_{T\in\cT_h} h_T$. The length of any edge $e\in \mathcal{E}_{h}$ will be denoted by $h_e$. Let $\cV_h$ denote the set of all the vertices of the triangles in $\cT_h$. The set of vertices on $\overline{\Gamma}_D,$ $\Gamma_N$ and $\Gamma_C$ are denoted by $\cV_h^D,$ $\cV_h^N$ and $\cV_h^C.$
Also, in the problem setting, we require the jump definitions of  scalar, vector valued functions and tensors. Let us define a broken Sobolev space
$$H^1(\Omega,\mathcal{T}_{h})=\{v\in L^2(\Omega): v|_{T}\in H^1(T) \;\;{\rm for~all~} T\in \mathcal{T}_h \}.$$ For any $e \in \mathcal{E}^{i}_{h}$, there are two triangles $T_{+}$ and $T_{-}$ such that $e=\partial T_{+} \cap \partial T_{-}$. Let
$\mathbf{n}_{+}$ be the unit normal of $e$ pointing from  $T_{+}$ to $T_{-}$ and let $\mathbf{n}_{-} =-\mathbf{n}_{+}$ (cf. Fig.\ref{Fig1}). For any $v \in H^1(\Omega,\mathcal{T}_{h})$, we define
the jump of $v$ on an edge $e$ by
$\sjump{v}=v_+\mathbf{n}_+ +v_-\mathbf{n}_-$
where $v_{\pm}=v|_{T_{\pm}}.$
\noindent
For $\mathbf{v} \in [H^1(\Omega,\mathcal{T}_{h})]^2$ we define the jump of $\mathbf{v}$ on $e \in \mathcal{E}^{i}_{h}.$ by
$\sjump{\mathbf{v}}=\mathbf{v}_+ \cdot \mathbf{n}_+ +\mathbf{v}_-\cdot \mathbf{n}_-.$
\noindent
 Similarly, for tensors $\tau  \in [H^1(\Omega,\mathcal{T}_{h})]^{2\times 2}$, the jump on $ e\in \mathcal{E}^i_{h}$ are defined by
$\sjump{\tau}=\tau_+\mathbf{n}_+ +\tau_-\mathbf{n}_-.$
\noindent
For notational convenience, we also define the jump on the boundary faces $e\in \mathcal{E}_{h}^b$ by modifying them appropriately. We use the definition of jump by understanding that $v_{-} = 0$ (similarly, $\mathbf{v}_{-} = 0$
and $\tau_{-} = 0$). 

Define the discrete space for velocity $\mathbf{V}_h\subset \mathbf{V}$ by
$$\mathbf{V}_h:=\{\mathbf{v}_h\in \mathbf{V}: \mathbf{v}_h|_T\in \mathbf{P}_1(T)\;\; \forall T \in \mathcal{T}_h\},$$
and the discrete space for pressure is
$$M_H:=\{p_H\in L^2(\Omega): p_H|_T\in P_0(T)\;\; \forall T \in \mathcal{T}_H\},$$
and the discrete control space $\mathbf{Q}_h\subset \mathbf{Q}$ by
$$\mathbf{Q}_h:=\{\mathbf{x}_h\in \mathbf{Q}: \mathbf{x}_h|_T\in \mathbf{P}_1(T),\;\; \forall T \in \mathcal{T}_h\},$$
where $\mathbf{P}_1(T)$ is the space of polynomials of degree less
than or equal to one on the triangle $T$. The discrete
admissible set of controls is defined by
$$\mathbf{Q}^h_{ad}:=\{\mathbf{x}_h\in \mathbf{Q}_h:\mathbf{y}_a\leq \mathbf{x}_h(z) \leq \mathbf{y}_b \text{ for all } z\in \cV_h^C\}.$$

It is easy to check that $\mathbf{Q}_{ad}^h\subset \mathbf{Q}_{ad}.$ Throughout the article, we assume that $C$ denotes a generic
positive constant that is independent of the mesh parameter $h$. A proof of the following proposition on the existence and uniqueness of the solution of discrete problem can be found in \cite[Proposition 3.1]{Gudi_Ramesh_CaMwA:2023}.	
\begin{proposition}[\bf Discrete Optimality System]\label{prop:Wellposedness-discrete_optmality}
There exists unique 
$((\mathbf{w}_h,p_H),(\bm{\phi}_h,r_H),\mathbf{y}_h)\in \big(\mathbf{V}_h\times M_H\big)\times \big(\mathbf{V}_h\times M_H\big) \times \mathbf{Q}^h_{ad}$ satisfying the following:
\begin{subequations}
\begin{align}
		\mathbf{u}_h&=\mathbf{w}_h+\mathbf{y}_h,\quad \mathbf{w}_h\in \mathbf{V}_h,\notag\\
		a(\mathbf{w}_h,\mathbf{z}_h)+b(\mathbf{z}_h,p_H) &=( \mathbf{f},\mathbf{z}_h)-a(\mathbf{y}_h,\mathbf{z}_h) \;\;\;
		{\rm for~all}\;\mathbf{z}_h \in \mathbf{V}_h,\label{eq:d_state_a}\\ 
		b(\mathbf{u}_h,q_H)&=0\;\quad {\rm for~all}\; q_H \in M_H ,\label{eq:d_state_b}\\
		a(\mathbf{z}_h,\bm{\phi}_h)-b(\mathbf{z}_h,r_H) &=( \mathbf{u}_h-\mathbf{u_d},\mathbf{z}_h) \;\;\;
		{\rm for~all}\;\mathbf{z}_h \in \mathbf{V}_h ,\label{eq:d_adjstate_a}\\ 
		b(\bm{\phi}_h,q_H)&=0\;\quad {\rm for~all}\; q_H \in M_H ,\label{eq:d_adjstate_b}\\
		\rho\, a(\mathbf{y}_h,\mathbf{x}_h-\mathbf{y}_h)\geq& a(\mathbf{x}_h-\mathbf{y}_h,\bm{\phi}_h)-b(\mathbf{x}_h-\mathbf{y}_h,r_H)\notag\\&-(\mathbf{u}_h-\mathbf{u_d},\mathbf{x}_h-\mathbf{y}_h)\;\quad {\rm for~all}\; \mathbf{x}_h\in
		\mathbf{Q}^h_{ad}.\label{eq:d_VI}
\end{align}
\end{subequations}
\end{proposition}

\section{{A~posteriori} Error Analysis}\label{Sto:DC_error_analysis}
This section is devoted to \textit{a~posteriori} error analysis.
Define reconstructions $\mathbf{R}\mathbf{w}\in \mathbf{V}$, $R_{0}p\in L^2(\Omega),$ $\bar{\mathbf{R}}\bm{\phi}\in \mathbf{V}$ and $\bar{R}_{0}r\in L^2(\Omega),$ $\mathbf{R}\mathbf{y}\in \mathbf{Q}_{ad}$  by
\begin{subequations}\label{apo1}
	\begin{align}
		\mathbf{R}\mathbf{u}=& \mathbf{R}\mathbf{w}+\mathbf{y}_h\nonumber\\
		a(\mathbf{R}\mathbf{w},\mathbf{v})+b(\mathbf{v},R_{0}p) = &{\langle \mathbf{f},\mathbf{z} \rangle}-a(\mathbf{y}_h,\mathbf{v}) \;~~\hspace{1.1 cm}{\rm for~all}\;\mathbf{v} \in \mathbf{V},\label{8}\\
		b(\mathbf{R}\mathbf{w},q)=&	-b(\mathbf{y}_h,q)\;~~\hspace{2.1 cm}{\rm for~all}\;q \in L^2(\Omega),\label{3.1b}\\
		a(\mathbf{z},\bar{\mathbf{R}}\bm{\phi})-b(\mathbf{z},\bar{R_{0}}r) =&{\langle \mathbf{u}_h-\mathbf{u}_d,\mathbf{z} \rangle}_{W} \;~~\hspace{1.6 cm}{\rm for~all}\;\mathbf{z} \in \mathbf{V},\label{9}\\
		b(\bar{\mathbf{R}}\bm{\phi},q)=&0\;~~\hspace{3.9 cm}{\rm for~all}\;q \in L^2(\Omega).\label{3.1d}\\
		\rho a(\mathbf{R}\mathbf{y},\mathbf{x}-\mathbf{R}\mathbf{y})\geq &  a(\mathbf{x}-\mathbf{R}\mathbf{y},\bm{\phi}_h)-b(\mathbf{x}-\mathbf{R}\mathbf{y},r_H)\label{3.1e}\\&- (\mathbf{u}_h-\mathbf{u}_d,\mathbf{x}-\mathbf{R}\mathbf{y})\;~~\hspace{0.5 cm}{\rm for~all}\;\mathbf{x} \in \mathbf{Q}_{ad}.\nonumber
	\end{align}
\end{subequations}

The well-posedness of the above system \eqref{apo1} follows from the facts that
the right-hand side of \eqref{8} is a bounded linear functional on $\mathbf{V}$,  the bilinear forms $a$ and $b$ are continuous, $a$ is elliptic and $b$ is inf-sup stable, and hence the system  \eqref{8}-\eqref{3.1b} has a unique solution  \cite[pp. 81]{Girault:1979} . Similarly, the system \eqref{9}-\eqref{3.1d} is well-posed.\\

\noindent
Subtraction of \eqref{apo1} from \eqref{conti_optimality_system} yields,
\begin{subequations}
	\begin{align}
		a(\mathbf{w}-\mathbf{R}\mathbf{w},\mathbf{z})+b(\mathbf{z},p-R_{0}p) &=-a(\mathbf{y}-\mathbf{y}_h,\mathbf{z})\; \;~~{\rm for~all}\;\mathbf{z} \in \mathbf{V},\label{recon_velocity} \\
		b(\mathbf{u}-\mathbf{R}\mathbf{u},q)&=0\;~~\hspace{2 cm}{\rm for~all}\;q \in L^2(\Omega),\label{3.2b}\\
		a(\mathbf{z},\bm{\phi}-\bar{\mathbf{R}}\bm{\phi})+b(\mathbf{z},\bar{R_{0}}r-r)& =(\mathbf{u}-\mathbf{u}_{h},\mathbf{z}) \;~~\hspace{0.3 cm}{\rm for~all}\;\mathbf{z} \in \mathbf{V},\label{11} \\
		b(\bm{\phi}-\bar{\mathbf{R}}\bm{\phi},q)& =0\;~~\hspace{2.4 cm}{\rm for~all}\;q \in L^2(\Omega).\label{3.2d}
	\end{align}	
\end{subequations}

\begin{theorem}[\bf Energy error estimate of control and $L^2$-estimate of velocity]\label{apost_energy_est_control}
	There holds,
	\begin{align*} 
		\rho\norm{\nabla(\mathbf{R}\mathbf{y}-\mathbf{y})}_{0,\Omega}+\norm{\mathbf{R}\mathbf{u}-\mathbf{u}}_{0,\Omega}\leq & C\big(\norm{\nabla(\mathbf{R}\mathbf{y}-\mathbf{y}_h)}_{0,\Omega}+\norm{\nabla(\mathbf{R}\mathbf{w}-\mathbf{w}_h)}_{0,\Omega}\\& + \norm{\nabla(\bm{\phi}_h-\bar{\mathbf{R}}\bm{\phi})}_{0,\Omega}+\norm{r_H-\bar{R}_0r}_{0,\Omega}\big).
	\end{align*}
\end{theorem}
\begin{proof} 
	Selecting $\mathbf{x=y}$ in \eqref{3.1e}, $\mathbf{x=Ry}$ in \eqref{eq:VI} and adding the result, we obtain
	\begin{align*}
		\rho a(\mathbf{R}\mathbf{y}-\mathbf{y},\mathbf{y}-\mathbf{R}\mathbf{y})\geq& a(\mathbf{y}-\mathbf{R}\mathbf{y},\bm{\phi}_h-\bm{\phi})-b(\mathbf{y}-\mathbf{R}\mathbf{y},r_H-r)-(\mathbf{u}_h-\mathbf{u},\mathbf{y}-\mathbf{R}\mathbf{y})\\
		\geq& a(\mathbf{y}-\mathbf{R}\mathbf{y},\bm{\phi}_h-\bar{\mathbf{R}}\bm{\phi})+a(\mathbf{y}-\mathbf{y}_h,\bar{\mathbf{R}}\bm{\phi}-\bm{\phi})+a(\mathbf{y}_h-\mathbf{R}\mathbf{y},\bar{\mathbf{R}}\bm{\phi}-\bm{\phi})\\&-b(\mathbf{y}-\mathbf{R}\mathbf{y},r_H-r)-(\mathbf{u}_h-\mathbf{u},\mathbf{y}-\mathbf{R}\mathbf{y})\\
		\geq& a(\mathbf{y}-\mathbf{R}\mathbf{y},\bm{\phi}_h-\bar{\mathbf{R}}\bm{\phi})+a(\mathbf{y}-\mathbf{y}_h,\bar{\mathbf{R}}\bm{\phi}-\bm{\phi})+a(\mathbf{y}_h-\mathbf{R}\mathbf{y},\bar{\mathbf{R}}\bm{\phi}-\bm{\phi})\\&-b(\mathbf{y}-\mathbf{R}\mathbf{y},r_H-r)-(\mathbf{u}_h-\mathbf{u},\mathbf{y}-\mathbf{R}\mathbf{y})\\
		\geq& a(\mathbf{y}-\mathbf{R}\mathbf{y},\bm{\phi}_h-\bar{\mathbf{R}}\bm{\phi})-b(\mathbf{y}-\mathbf{y}_h,r-\bar{R}_0r)+a(\mathbf{y}_h-\mathbf{R}\mathbf{y},\bar{\mathbf{R}}\bm{\phi}-\bm{\phi})\\&-b(\mathbf{y}-\mathbf{R}\mathbf{y},r_H-r)-(\mathbf{u}-\mathbf{u}_h,\mathbf{y}+\mathbf{w}-\mathbf{R}\mathbf{y}-\mathbf{R}\mathbf{w})\\
		\geq& a(\mathbf{y}-\mathbf{R}\mathbf{y},\bm{\phi}_h-\bar{\mathbf{R}}\bm{\phi})-b(\mathbf{y}-\mathbf{y}_h,r-\bar{R}_0r)+a(\mathbf{y}_h-\mathbf{R}\mathbf{y},\bar{\mathbf{R}}\bm{\phi}-\bm{\phi})\\&-b(\mathbf{y}-\mathbf{R}\mathbf{y},r_H-r)+\norm{\mathbf{u}-\mathbf{R}\mathbf{u}}_{0,\Omega}^2+(\mathbf{R}\mathbf{u}-\mathbf{u}_h,\mathbf{u}-\mathbf{R}\mathbf{u})\\&+(\mathbf{u}-\mathbf{u}_h,\mathbf{y}_h-\mathbf{R}\mathbf{y}).
	\end{align*}
	Now we have
	\begin{align*}
		\rho\norm{\nabla(\mathbf{R}\mathbf{y}-\mathbf{y})}_{0,\Omega}^2+\norm{\mathbf{u}-\mathbf{R}\mathbf{u}}_{0,\Omega}^2\leq& -a(\mathbf{y}-\mathbf{R}\mathbf{y},\bm{\phi}_h-\bar{\mathbf{R}}\bm{\phi})+b(\mathbf{y}-\mathbf{y}_h,r-\bar{R}_0r)\\&-a(\mathbf{y}_h-\mathbf{R}\mathbf{y},\bar{\mathbf{R}}\bm{\phi}-\bm{\phi})+b(\mathbf{y}-\mathbf{R}\mathbf{y},r_H-r)\\&-(\mathbf{R}\mathbf{u}-\mathbf{u}_h,\mathbf{u}-\mathbf{R}\mathbf{u})-(\mathbf{u}-\mathbf{u}_h,\mathbf{y}_h-\mathbf{R}\mathbf{y})\\
		\leq& -a(\mathbf{y}-\mathbf{R}\mathbf{y},\bm{\phi}_h-\bar{\mathbf{R}}\bm{\phi})+b(\mathbf{R}\mathbf{y}-\mathbf{y}_h,r-\bar{R}_0r)\\&+b(\mathbf{y}-\mathbf{R}\mathbf{y},r_H-\bar{R}_0r)-a(\mathbf{y}_h-\mathbf{R}\mathbf{y},\bar{\mathbf{R}}\bm{\phi}-\bm{\phi})\\&-(\mathbf{R}\mathbf{u}-\mathbf{u}_h,\mathbf{u}-\mathbf{R}\mathbf{u})-(\mathbf{R}\mathbf{u}-\mathbf{u}_h,\mathbf{y}_h-\mathbf{R}\mathbf{y})\\&+(\mathbf{u}-\mathbf{R}\mathbf{u},\mathbf{y}_h-\mathbf{R}\mathbf{y}).
	\end{align*}
	Applying Cauchy-Schwarz inequality and Young's inequality we arrive at the desired estimate. 
\end{proof}
\begin{theorem}[\bf Energy error estimate of velocity]\label{energy_est_state}
	There holds,
	\begin{align*}
		\norm{\nabla(\mathbf{u}-\mathbf{R}\mathbf{u})}_{0,\Omega}\leq& C\big(\norm{\nabla(\mathbf{R}\mathbf{y}-\mathbf{y}_h)}_{0,\Omega}+\norm{\nabla(\mathbf{R}\mathbf{w}-\mathbf{w}_h)}_{0,\Omega}+\norm{\nabla(\bm{\phi}_h-\bar{\mathbf{R}}\bm{\phi})}_{0,\Omega}\\&+\norm{r_H-\bar{R}_0r}_{0,\Omega}\big).
	\end{align*}	
\end{theorem}

\begin{proof}
	The splitting  $\mathbf{u=w+y}$ and $\mathbf{Ru=Rw+y_h}$ yields   	
	\begin{align}\label{Est2}
		\norm{\nabla(\mathbf{u}-\mathbf{R}\mathbf{u})}_{0,\Omega}\leq \norm{\nabla(\mathbf{w}-\mathbf{R}\mathbf{w})}_{0,\Omega}+\norm{\nabla(\mathbf{y}-\mathbf{y}_h)}_{0,\Omega}.
	\end{align}
	Now we need to estimate the term $\norm{\nabla(\mathbf{w}-\mathbf{R}\mathbf{w})}_{0,\Omega}.$ A selection $\mathbf{z}=\mathbf{w}-\mathbf{R}\mathbf{w}$ in \eqref{recon_velocity} provides
	\begin{align}\label{recon_aux}
		\norm{\nabla(\mathbf{w}-\mathbf{R}\mathbf{w})}_{0,\Omega}^2+b(\mathbf{w}-\mathbf{R}\mathbf{w},p-R_0p)=-a(\mathbf{y}-\mathbf{y}_h,\mathbf{w}-\mathbf{R}\mathbf{w}).
	\end{align}
	Using the fact that $b(\mathbf{w}-\mathbf{R}\mathbf{w},p-R_0p)=-b(\mathbf{y}-\mathbf{y}_h,p-R_0p)$ and appplying Cauchy-Schwarz and Young's inequality in \eqref{recon_aux}, we obtain
	\begin{align}\label{Est}
		(1-\frac{\epsilon}{2})\norm{\nabla(\mathbf{w}-\mathbf{R}\mathbf{w})}_{0,\Omega}^2\leq  \frac{1}{\epsilon}\norm{\nabla(\mathbf{y}-\mathbf{y}_h)}_{0,\Omega}^2+\frac{\epsilon}{2}\norm{p-R_0p}_{0,\Omega}^2.
	\end{align}
	The estimate of $\norm{p-R_0p}_{0,\Omega}$ is in the following:
	\begin{align}\label{Est1}
		\norm{p-R_0p}_{0,\Omega}&\leq \sup_{\mathbf{v}\in \mathbf{V}}\frac{b(\mathbf{v},p-R_0p)}{\norm{v}_\mathbf{V}}\nonumber\\
		&\leq \sup_{\mathbf{v}\in \mathbf{V}}\frac{-a(\mathbf{y}-\mathbf{y_h},\mathbf{v})-a(\mathbf{w}-\mathbf{R}\mathbf{w},\mathbf{v})}{\norm{v}_\mathbf{V}}\nonumber\\
		&\leq\norm{\nabla(\mathbf{y}-\mathbf{y}_h)}_{0,\Omega}+\norm{\nabla(\mathbf{w}-\mathbf{R}\mathbf{w})}_{0,\Omega}.
	\end{align}
	Using \eqref{Est1} in \eqref{Est}, we obtain 
	$\norm{\nabla(\mathbf{w}-\mathbf{R}\mathbf{w})}_{0,\Omega}\leq C \norm{\nabla(\mathbf{y}-\mathbf{y}_h)}_{0,\Omega}.$
	 Hence from \eqref{Est2} we get
	\begin{align}\label{Est3}
		\norm{\nabla(\mathbf{u}-\mathbf{R}\mathbf{u})}_{0,\Omega}&\leq C \norm{\nabla(\mathbf{y}-\mathbf{y}_h)}_{0,\Omega}
	\end{align} 
	Introducing the reconstruction $\mathbf{R}\mathbf{y}$ in \eqref{Est3} we obtain
	\begin{align}\label{es1}
		\norm{\nabla(\mathbf{u}-\mathbf{R}\mathbf{u})}_{0,\Omega}&\leq C \big(\norm{\nabla(\mathbf{y}-\mathbf{R}\mathbf{y})}_{0,\Omega}+\norm{\nabla(\mathbf{R}\mathbf{y}-\mathbf{y}_h)}_{0,\Omega}\big).
	\end{align}
	Substituting the estimate of $\norm{\nabla(\mathbf{R}\mathbf{y}-\mathbf{y})}_{0,\Omega}$(from Theorem \ref{apost_energy_est_control}) in \eqref{es1}, we obtain the desired estimate.
\end{proof}
\begin{theorem}[\bf Energy error estimate of adjoint velocity]
	There holds,
	\begin{align*}
		\norm{\nabla(\bm{\phi}-\mathbf{R}\bm{\phi})}_{0,\Omega}\leq& C\big(\norm{\nabla(\mathbf{R}\mathbf{y}-\mathbf{y}_h)}_{0,\Omega}+\norm{\nabla(\mathbf{R}\mathbf{w}-\mathbf{w}_h)}_{0,\Omega}+\norm{\nabla(\bm{\phi}_h-\bar{\mathbf{R}}\bm{\phi})}_{0,\Omega}\\&+\norm{r_H-\bar{R}_0r}_{0,\Omega}\big).
	\end{align*}	
\end{theorem}
\begin{proof}
	Selecting $\mathbf{z}=\bm{\phi}-\mathbf{R}\bm{\phi}$ in \eqref{11} we have the following:
	\begin{align}\label{est4}
		\norm{\nabla(\bm{\phi}-\bar{\mathbf{R}}\bm{\phi})}_{0,\Omega}^2-b(\bm{\phi}-\bar{\mathbf{R}}\bm{\phi},r-\bar{R}_0r)=(\mathbf{u}-\mathbf{u}_h,\bm{\phi}-\bar{\mathbf{R}}\bm{\phi}).
	\end{align}
	Using the fact that
	$b(\bm{\phi}-\bar{\mathbf{R}}\bm{\phi},r-\bar{R}_0r)=0$ and applying Cauchy-Schwarz inequality in \eqref{est4}, we obtain 
	\begin{align}
		\norm{\nabla(\bm{\phi}-\bar{\mathbf{R}}\bm{\phi})}_{0,\Omega}&\leq \norm{\mathbf{u}-\mathbf{u}_h}_{0,\Omega}\notag\\
		&\leq\norm{\mathbf{u}-\mathbf{Ru}}_{0,\Omega}+\norm{\mathbf{Ru}-\mathbf{u}_h}_{0,\Omega}\notag\\
		&\leq\norm{\mathbf{u}-\mathbf{Ru}}_{0,\Omega}+\norm{\mathbf{Rw}-\mathbf{w}_h}_{0,\Omega}\notag\\
		&\leq\norm{\mathbf{u}-\mathbf{Ru}}_{0,\Omega}+\norm{\nabla(\mathbf{Rw}-\mathbf{w}_h)}_{0,\Omega}.\label{apost_adj_estm}
	\end{align} 
	In the above, we have used the fact that $\mathbf{u}_h=\mathbf{w}_h+\mathbf{y}_h$, $\mathbf{Ru}=\mathbf{Rw}+\mathbf{y}_h$ and Poincar\'e inequality. Using the estimates of $\norm{\mathbf{u}-\mathbf{Ru}}_{0,\Omega}$ from Theorem \ref{apost_energy_est_control} in \eqref{apost_adj_estm} we obtain the desired result.
\end{proof}
\noindent
In the next theorem we state the estimate of pressure and adjoint pressure the proof can be easily derived by using the inf-sup condition.
\begin{theorem}[\bf Error estimate of pressure and adjoint pressure]\label{PressureEstm}
	There holds,
	\begin{align*}
		\norm{p-R_0p}_{0,\Omega}\leq& C\big(\norm{\nabla(\mathbf{R}\mathbf{y}-\mathbf{y}_h)}_{0,\Omega}+\norm{\nabla(\mathbf{R}\mathbf{w}-\mathbf{w}_h)}_{0,\Omega}+\norm{\nabla(\bm{\phi}_h-\bar{\mathbf{R}}\bm{\phi})}_{0,\Omega}\\&+\norm{r_H-\bar{R}_0r}_{0,\Omega}\big),
	\end{align*}
and
\begin{align*}
	\norm{r-\bar{R}_0r}_{0,\Omega}\leq& C\big(\norm{\nabla(\mathbf{R}\mathbf{y}-\mathbf{y}_h)}_{0,\Omega}+\norm{\nabla(\mathbf{R}\mathbf{w}-\mathbf{w}_h)}_{0,\Omega}+\norm{\nabla(\bm{\phi}_h-\bar{\mathbf{R}}\bm{\phi})}_{0,\Omega}\\&+\norm{r_H-\bar{R}_0r}_{0,\Omega}\big).
	\end{align*}	
\end{theorem}

\noindent
Combining all the above results, Theorem \ref{apost_energy_est_control}-\ref{PressureEstm}, we obtain the following theorem:
\begin{theorem}\label{apost_aux}
	There holds
	\begin{align*}
		\norm{\nabla(\mathbf{u}-\mathbf{u}_h)}_{0,\Omega}+&\norm{p-p_H}_{0,\Omega}+\norm{\nabla(\bm{\phi}-\bm{\phi}_h)}_{0,\Omega}+\norm{r-r_H}_{0,\Omega}+\norm{\nabla(\mathbf{y}-\mathbf{y}_h)}_{0,\Omega}\leq C \big(\norm{\nabla(\mathbf{R}\mathbf{y}-\mathbf{y}_h)}_{0,\Omega}\\&+\norm{\nabla(\mathbf{R}\mathbf{w}-\mathbf{w}_h)}_{0,\Omega}+\norm{p_H-R_0p}_{0,\Omega}+\norm{\nabla(\bm{\phi}_h-\bar{\mathbf{R}}\bm{\phi})}_{0,\Omega}+\norm{r_H-\bar{R}_0r}_{0,\Omega}\big).
	\end{align*}
\end{theorem}

Before going to derive the \textit{a posteriori} error estimator we need to define some preliminary definitions which are given in the following:
Let $p \in \cV_h$ be a node and $\omega_p$ be the node patch of $p$ and define $h_p=\text{diam}\; \omega_p.$ Denote $\gamma_p$ be the union of all sides of $\bar{\omega}_p$ and union of interior sides of $\bar{\omega}_p$ are denoted by $\gamma_{p,I}.$ Given any $p\in \cV_h^C$ we subdivide the intersections between $\partial\Omega$ and $\partial\omega_p$ in the three following sets: 
\begin{align*}
	\gamma_{p,C}&:=\Gamma_C \cap \partial \omega_p\\
	\gamma_{p,D}&:=\Gamma_D \cap \partial \omega_p\\
	\gamma_{p,N}&:=\Gamma_N \cap \partial \omega_p.
\end{align*}

 We define the discrete contact stress by $\hat{\bm{\sigma}}(\mathbf{y}_h):=\rho\frac{\partial \mathbf{y}_h}{\partial n}-\frac{\partial \bm{\phi}_h}{\partial
	n}-r_H\mathbf{n},$ clearly it is a vector quantity so the sign of this quantity would be in the sense of componentwise. We classify the nodes on $\Gamma_C$ as follows.

We classify the actual contact nodes $p\in\cV_h^C$
with $\mathbf{y}_h(p)=\mathbf{y}_a(p)(=\mathbf{y}_a)$ in two different categories. At so-called full-contact nodes $p\in \cV^{fC}_{h,a}$ if the discrete solution satisfies $\mathbf{y}_h=\mathbf{y}_a,$ $\hat{\bm{\sigma}}(\mathbf{y}_h):=(\hat{\sigma}_1(\mathbf{y}_h),\hat{\sigma}_2(\mathbf{y}_h))\geq \mathbf{0}$ on $\gamma_{p,C},$ which means that the conditions of actual contact are satisfied. The
remaining actual contact nodes with $\mathbf{y}_h(p)=\mathbf{y}_a(p)$ are called semi-contact nodes and the set is denoted by $\cV^{sC}_{h,a}.$

Similarly, we classify the actual contact nodes $p\in\cV_h^C$
with $\mathbf{y}_h(p)=\mathbf{y}_b(p)(=\mathbf{y}_b)$ in two different categories. At so-called full-contact nodes $p\in \cV^{fC}_{h,b}$ if the discrete solution satisfies $\mathbf{y}_h=\mathbf{y}_b,$ $\hat{\bm{\sigma}}(\mathbf{y}_h)\leq \mathbf{0}$ on $\gamma_{p,C},$ which means that the conditions of actual contact are satisfied. The
remaining actual contact nodes with $\mathbf{y}_h(p)=\mathbf{y}_b(p)$ are called semi-contact nodes and the set is denoted by $\cV^{sC}_{h,b}.$

\noindent
Also, we define $\cV_{*}^C=\{p\in \cV ^C_h: \mathbf{y}_a<\mathbf{y}_h(p)<\mathbf{y}_b\}.$  
It is clear that $$\cV^C_h=\cV^{fC}_{h,a}\cup\cV^{fC}_{h,b}\cup\cV^{sC}_{h,a}\cup\cV^{sC}_{h,b}\cup\cV_{*}^C.$$ 
We define the residual $\mathbf{R}_h$ by the following:
\begin{align*}
	(\mathbf{R}_h,\mathbf{x})_{-1,1}=a(\mathbf{x},\bm{\phi}_h)-b(\mathbf{x},r_H)-(\mathbf{u}_h-\mathbf{u}_d,\mathbf{x})-\rho a(\mathbf{y}_h,\mathbf{x})\quad \text{for~all}\; \mathbf{x}\in \mathbf{Q}.
\end{align*}
Also, for the \textit{a posteriori} estimator we need to define the the Lagrange multiplier $\bm{\lambda}$ by the following:
\begin{align*}
	(\bm{\lambda},\mathbf{x})_{-1,1}=a(\mathbf{x},\bm{\phi}_h)-b(\mathbf{x},r_H)-(\mathbf{u}_h-\mathbf{u}_d,\mathbf{x})-\rho a(\mathbf{Ry},\mathbf{x})\quad \text{for~all}\; \mathbf{x}\in \mathbf{Q}.
\end{align*}
Here, $(\cdot,\cdot)_{-1,1}$ denotes the duality pairing between $\mathbf{Q}$ and its dual. Clearly, $(\bm{\lambda},\mathbf{x}-\mathbf{Ry})_{-1,1}\leq 0$ for all $\mathbf{x}\in \mathbf{Q}_{ad}.$
Define the discrete Lagrange multiplier $\bm{\lambda}_h$ by the following:
\begin{align*}
	(\bm{\lambda}_h,\mathbf{x}_h)_{-1,1}=a(\mathbf{x}_h,\bm{\phi}_h)-b(\mathbf{x}_h,r_H)-(\mathbf{u}_h-\mathbf{u}_d,\mathbf{x}_h)-\rho a(\mathbf{y}_h,\mathbf{x}_h)\quad \text{for~all}\; \mathbf{x}_h\in \mathbf{Q}_h.
\end{align*}
It is clear that, $(\bm{\lambda}_h,\mathbf{x}_h)_{-1,1}=(\mathbf{R}_h,\mathbf{x}_h)_{-1,1}$ for all $\mathbf{x}_h\in \mathbf{Q}^h_{ad}$ and
\begin{align*}
	(\bm{\lambda}_h,\mathbf{x}_h-\mathbf{y}_h)_{-1,1}\leq 0 \quad \text{for~all}\; \mathbf{x}_h\in \mathbf{Q}^h_{ad}.
\end{align*}
In order to investigate $\bm{\lambda}_h$ further, we use the partition of unity and integration by parts
		\begin{align*}
			(\bm{\lambda}_h,\mathbf{x}_h)_{-1,1}=&\sum_{i=1}^{2}\sum_{p\in \cV_{h}}(\bm{\lambda}_h,x_{h,i}(p)\phi_p\mathbf{e}_i)_{-1,1}\\
		=&\sum_{i=1}^{2}\sum_{p\in \cV_{h}} \int_{\omega_p}(-\Delta \bm{\phi}_h-\nabla r_H-(\mathbf{u}_h-\mathbf{u}_d)+\rho \Delta \mathbf{y}_h)\cdot x_{h,i}(p)\phi_p\mathbf{e}_i\\&+\sum_{i=1}^{2}\sum_{p\in \cV_{h}}\int_{\gamma_{p,I}}\sjump{\nabla\bm{\phi}_h+r_H\mathbb{I}-\rho \nabla\mathbf{y}_h}\cdot x_{h,i}(p)\phi_p\mathbf{e}_i\\&-\sum_{i=1}^{2}\sum_{p\in \cV^C_{h}} \int_{\gamma_{p,C}} \hat{\bm{\sigma}}(\mathbf{y}_h)\cdot x_{h,i}(p)\phi_p\mathbf{e}_i,  
	\end{align*}
where $\{\mathbf{e}_1,\mathbf{e}_2\}$ standard basis of $\mathbb{R}^2,$ and $\phi_p$ is the Lagrange basis at node $p$. This motivates the representation
$(\bm{\lambda}_h,\mathbf{x}_h)_{-1,1}=\sum_{i=1}^{2}(\lambda_{h,i},x_{h,i})_{-1,1}.$

Now for the \textit{a posteriori} estimator we replace the residual $\mathbf{R}_h$ by a Galerkin functional, whose abstract definition is given by
\begin{align}\label{galerkin_residual}
	(\mathbf{G}_h,\mathbf{x})_{-1,1}:=\rho a(\mathbf{Ry}-\mathbf{y}_h,\mathbf{x})+(\bm{\lambda}-\tilde{\bm{\lambda}}_h,\mathbf{x})_{-1,1}
\end{align} 
for all $\mathbf{x}\in \mathbf{Q},$ where $\tilde{\bm{\lambda}}_h$ is an approximation of $\bm{\lambda}_h$ and it depends on the discrete solution, data and reflects the properties of $\bm{\lambda},$ details can be found in \cite{Fierro_Veeser2003}. We call it quasi-discrete contact force
density.

\noindent
By using the partition of unity we have the following definition of $\tilde{\bm{\lambda}}_h,$ for all $\bm{\psi}\in \mathbf{Q}$:
\begin{align}\label{Def_lambda_tilde}
	(\tilde{\bm{\lambda}}_h,\bm{\psi})_{-1,1}:=\sum_{i=1}^{2}(\tilde{\lambda}_{h,i},\psi_i)_{-1,1}=\sum_{i=1}^{2}\sum_{p\in \cV_h^C}(\tilde{\lambda}^{p}_{h,i},\psi_i\phi_p)_{-1,1}
\end{align}
and it adjusts the local contributions so that, on one hand, the Galerkin functional is
prepared for the derivation of an upper bound, and on the other hand, tries to maximize
the cancellation within $\mathbf{G}_h.$  

\noindent
For semi-contact nodes $p\in \cV^{sC}_h:=\cV^{sC}_{h,a}\cup \cV^{sC}_{h,b},$ we define
\begin{align*}
	(\tilde{\lambda}^{p}_{h,i},\psi_i\phi_p)_{-1,1}:=(R_{h,i},\phi_p)_{-1,1}m_p(\psi_i)=s_{p,i} m_p(\psi_i)\int_{\gamma_{p,C}}\phi_p,
\end{align*}
where $s_{p,i}:=\frac{(\lambda_{h,i},\phi_p)_{-1,1}}{\int_{\gamma_{p,C}}\phi_p},$ $i=1,2.$ The constant $s_{p,i}$ is the nodal value of the discrete contact force density obtained by lumping the boundary mass matrix and $m_p(\psi_i)$ is defined below.
Sign of $s_{p,i}$ is the following:\\
i) For $p\in \Gamma_C$ with $y_{h,i}(p)=y^{i}_a$ then $s_{p,i}\leq 0.$\\
ii) For $p\in \Gamma_C$ with $y_{h,i}(p)=y^{i}_b$ then $s_{p,i}\geq 0.$\\
iii) For $p\in \cV_{*}^C$ then $s_{p,i}=0.$

\noindent
For full-contact nodes $p\in \cV^{fC}_h:= \cV^{fC}_{h,a}\cup\cV^{fC}_{h,b},$ we define 
\begin{align*}
	(\tilde{\lambda}^{p}_{h,i},\psi_i\phi_p)_{-1,1}:=&(\tilde{R}_{h,i},\phi_p)_{-1,1}m_p(\psi_i)-\int_{\gamma_{p,C}}\hat{\sigma}_i(\mathbf{y}_h)\psi_i\phi_p\nonumber\\=&s_{p,i} m_p(\psi_i)\int_{\gamma_{p,C}}\phi_p-\int_{\gamma_{p,C}}\hat{\sigma}_i(\mathbf{y}_h)(\psi_i-m_p(\psi_i))\phi_p,		
\end{align*}
where,
\begin{align*}
	(\tilde{R}_{h,i},\theta)_{-1,1}:=(R_{h,i},\theta)_{-1,1}+\int_{\Gamma_{C}}\hat{\sigma}_i(\mathbf{y}_h)\theta.
\end{align*}
We need to specify the choices of $m_p(\psi_i)$ for semi- and full-contact nodes.

\noindent
For full-contact nodes $p\in\cV^{fC}_{h,b}$ we use 
\begin{align}\label{c_p_def_full_contact_b}
	m_p(\psi_i)=\max_{e\subset \gamma_{p,C}} \frac{\int_{e}\psi_i\phi_p}{\int_{e}\phi_p}
\end{align}
 and for full-contact nodes $p\in \cV^{fC}_{h,a}$ we use 
\begin{align}\label{c_p_def_full_contact_a}
	m_p(\psi_i)=\min_{e\subset \gamma_{p,C}} \frac{\int_{e}\psi_i\phi_p}{\int_{e}\phi_p}.
\end{align}
This choice is important for the derivation of the upper bound, see for instance \eqref{fourth_term} and \eqref{sixth_term}.

For semi-contact nodes $p\in \cV^{sC}_{h,a}\cup\cV^{sC}_{h,b},$ we set 
\begin{align}\label{c_p_def_semi_contact}
	m_p(\psi_i)=\frac{\int_{\tilde{\gamma}_{p,C}}\psi_i\phi_p}{\int_{\tilde{\gamma}_{p,C}}\phi_p},
\end{align}
where $\tilde{\gamma}_{p,C}$ is a strict subset of $\gamma_{p,C}.$ 
Now we are ready to derive the \textit{a posteriori} error estimators.

\begin{theorem}[\bf \textit{A~posteriori} error estimator]\label{aposteriori error estimator} It holds,
	\begin{align*}
		\lVert \nabla(\mathbf{y}-\mathbf{y}_h)\rVert_{0,\Omega}+&\lVert \nabla(\mathbf{u}-\mathbf{u}_{h})\rVert_{0,\Omega} +\lVert p-p_H\rVert_{0,\Omega}+\lVert \nabla (\bm{\phi}-\bm{\phi}_{h})\rVert_{0,\Omega}\\
		&+\lVert r-r_H\rVert_{0,\Omega}+\norm{\bm{\lambda}-\tilde{\bm{\lambda}}_h}_{-1,\Omega}\lesssim \big( \eta_{(\mathbf{u},p)}+\eta_{(\bm{\phi},r)}+\eta_{\mathbf{y}}\big),
	\end{align*}
	where the estimators are defined as
	\begin{align*}
		\eta_{(\mathbf{u},p)}^2=&\sum_{T\in \mathcal{T}_{h}} h_{T}^2 \lVert \mathbf{f}+\Delta \mathbf{u}_h-\nabla p_H\rVert^2_{0,T}+\sum_{e\in \mathcal{E}_{h}^{i}} \lVert h_e^{\frac{1}{2}}\sjump{p_{H}\I-\nabla{\mathbf{u}_{h}}}\rVert_{0,e}^2\\&+\sum_{e\in \mathcal{E}_{h}^{b,N}} \lVert h_e^{\frac{1}{2}}\sjump{p_{H}\I-\nabla{\mathbf{u}_{h}}}\rVert_{0,e}^2+ \sum_{T\in \mathcal{T}_h} \norm{\nabla \cdot \mathbf{u}_h}_{0,T}^2
	\end{align*}
	and
	\begin{align*}
		\eta_{(\bm{\phi},r)}^2=&\sum_{T\in \mathcal{T}_{h}} h_{T}^2 \lVert \Delta \bm{\phi}_h+\nabla r_H+\mathbf{u}_h-\mathbf{u}_d\rVert^2_{0,T}+\sum_{e\in \mathcal{E}_{h}^{i}} \lVert h_e^{\frac{1}{2}}\sjump{r_{H}\I+\nabla{\bm{\phi}_{h}}}\rVert_{0,e}^2\\&+\sum_{e\in \mathcal{E}_{h}^{b,N}} \lVert h_e^{\frac{1}{2}}\sjump{r_{H}\I+\nabla{\bm{\phi}_{h}}}\rVert_{0,e}^2+ \sum_{T\in \mathcal{T}_h} \norm{\nabla \cdot \bm{\phi}_h}_{0,T}^2
	\end{align*}

\begin{align*}
	\eta_{\mathbf{y}}^2=&
	\sum_{p\in \cV_{h}} \big(\underbrace{h_p^2\norm{\bm{\mathcal{R}}(\mathbf{y}_h)}^2_{0,\omega_p}}_{\eta_\mathbf{y}^a}+\underbrace{h_p \norm{\bm{\mathcal{J}}^{I}(\mathbf{y}_h)}_{0,\gamma_{p,I}}^2}_{\eta_\mathbf{y}^b}\big)+\sum_{p\in \cV^C_{h}\setminus \cV^{fC}_{h}}\underbrace{h_p \norm{\hat{\bm{\sigma}}(\mathbf{y}_h)}_{0,\gamma_{p,C}}^2}_{\eta_\mathbf{y}^c}\nonumber\\&+\sum_{i=1}^{2}\big[\sum_{p\in \cV^{sC}_{h,a}}\underbrace{ s_{p,i} d^{a}_{p,i}}_{\eta_\mathbf{y}^d}+ \sum_{p\in \cV^{sC}_{h,b}} \underbrace{s_{p,i} d^{b}_{p,i}}_{\eta_\mathbf{y}^e}\big],
\end{align*}

where, $\bm{\mathcal{R}}(\mathbf{y}_h):=(\mathcal{R}_1(\mathbf{y}_h),\mathcal{R}_2(\mathbf{y}_h))=-\Delta \bm{\phi}_h-\nabla r_H-(\mathbf{u}_h-\mathbf{u}_d)+\rho \Delta \mathbf{y}_h,$  $\bm{\mathcal{J}}^{I}(\mathbf{y}_h):=(\mathcal{J}^{I}_1(\mathbf{y}_h),\mathcal{J}^{I}_2(\mathbf{y}_h))=\sjump{\nabla\bm{\phi}_h+r_H\I-\rho \nabla\mathbf{y}_h},$  $d^{a}_{p,i}=\int_{\tilde{\gamma}_{p,C}}(y^i_a-y_{h,i})\phi_p.$ and $d^{b}_{p,i}=\int_{\tilde{\gamma}_{p,C}}(y^i_b-y_{h,i})\phi_p,$ $\I$ be the $2\times2$ identity matrix. and the set $\tilde{\gamma}_{p,C}$ is a strict subset of $\gamma_{p,C}.$

\end{theorem}
\begin{proof} From Theorem \ref{apost_aux}, we have
	\begin{align}\label{apostaux}
		\lVert \nabla(\mathbf{y}-\mathbf{y}_h)\rVert_{0,\Omega}+&\lVert \nabla(\mathbf{u}-\mathbf{u}_{h})\rVert_{0,\Omega} +\lVert p-p_H\rVert_{0,\Omega}+\lVert \nabla (\bm{\phi}-\bm{\phi}_{h})\rVert_{0,\Omega}+\norm{\bm{\lambda}-\tilde{\bm{\lambda}}_h}_{-1,\Omega}\nonumber\\&+\lVert r-r_H\rVert_{0,\Omega}\lesssim \norm{\nabla(\mathbf{R}\mathbf{y}-\mathbf{y}_h)}_{0,\Omega}+\norm{\bm{\lambda}-\tilde{\bm{\lambda}}_h}_{-1,\Omega}+\norm{\nabla(\mathbf{R}\mathbf{w}-\mathbf{w}_h)}_{0,\Omega}\nonumber\\&+\norm{p_H-R_0p}_{0,\Omega}+\norm{\nabla(\bar{\mathbf{R}}\bm{\phi}-\bm{\phi}_h)}_{0,\Omega}+\norm{r_H-\bar{R}_0r}_{0,\Omega}.
	\end{align}
	The {\it a~posteriori} error analysis in  \cite[Theorem 3.1]{HSW:2005:Err.Stokes}  and \cite[Section 7]{CarstensenEigelLoebhardetal.2010}  gives the following error estimators of state and adjoint state variables:	
	\begin{align}\label{apost_state}
		\norm{\nabla(\mathbf{R}\mathbf{u}-\mathbf{u}_h)}^2_{0,\Omega}+\norm{p_H-R_0p}^2_{0,\Omega}\lesssim &\sum_{T\in \mathcal{T}_{h}} h_{T}^2 \lVert \mathbf{f}+\Delta \mathbf{u}_h-\nabla p_H\rVert^2_{0,T}+\sum_{e\in \mathcal{E}_{h}^{i}} \lVert h_e^{\frac{1}{2}}\sjump{p_{H}\I-\nabla{\mathbf{u}_{h}}}\rVert_{0,e}^2\nonumber\\&+\sum_{e\in \mathcal{E}_{h}^{b,N}} \lVert h_e^{\frac{1}{2}}\sjump{p_{H}\I-\nabla{\mathbf{u}_{h}}}\rVert_{0,e}^2+ \sum_{T\in \mathcal{T}_h} \norm{\nabla \cdot \mathbf{u}_h}_{0,T}^2
	\end{align}
	\begin{align}\label{apost_adj_state}
		\norm{\nabla(\bar{\mathbf{R}}\bm{\phi}-\bm{\phi}_h)}^2_{0,\Omega}+\norm{r_H-\bar{R}_0r}^2_{0,\Omega}\lesssim & \sum_{T\in \mathcal{T}_{h}} h_{T}^2 \lVert \Delta \bm{\phi}_h+\nabla r_H+\mathbf{u}_h-\mathbf{u}_d\rVert^2_{0,T}+\sum_{e\in \mathcal{E}_{h}^{i}} \lVert h_e^{\frac{1}{2}}\sjump{r_{H}\I+\nabla{\bm{\phi}_{h}}}\rVert_{0,e}^2\nonumber\\&+\sum_{e\in \mathcal{E}_{h}^{b,N}} \lVert h_e^{\frac{1}{2}}\sjump{r_{H}\I+\nabla{\bm{\phi}_{h}}}\rVert_{0,e}^2+ \sum_{T\in \mathcal{T}_h} \norm{\nabla \cdot \bm{\phi}_h}_{0,T}^2.
	\end{align}
	We denote the right hand side of \eqref{apost_state} by $\eta^2_{(\mathbf{u},p)}$ and the right hand side of \eqref{apost_adj_state} by $\eta^2_{(\bm{\phi},r)}.$ The splitting $\mathbf{Ru}=\mathbf{Rw}+\mathbf{y}_h$ and $\mathbf{u}_h = \mathbf{w}_h+\mathbf{y}_h,$ yields $\norm{\nabla(\mathbf{R}\mathbf{u}-\mathbf{u}_h)}_{0,\Omega}=\norm{\nabla(\mathbf{R}\mathbf{w}-\mathbf{w}_h)}_{0,\Omega}.$ Thus we have estimates for all the last four terms of \eqref{apostaux}. Therefore we only need to estimate $\norm{\nabla(\mathbf{R}\mathbf{y}-\mathbf{y}_h)}_{0,\Omega}+\norm{\bm{\lambda}-\tilde{\bm{\lambda}}_h}_{-1,\Omega}.$ Estimator for the control variable needs some special care.
	\noindent
	To find an upper bound of the term $\norm{\nabla(\mathbf{Ry}-\mathbf{y}_h)}_{0,\Omega},$ 
	we choose $\mathbf{x}=\mathbf{Ry}-\mathbf{y}_h$ in 
	\eqref{galerkin_residual}
	\begin{align}
		\rho \norm{\nabla(\mathbf{Ry}-\mathbf{y}_h)}_{0,\Omega}^2&= \rho \; a(\mathbf{Ry}-\mathbf{y}_h,\mathbf{Ry}-\mathbf{y}_h)\notag\\
		&= \mathbf{G}_h(\mathbf{Ry}-\mathbf{y}_h)-(\bm{\lambda}-\tilde{\bm{\lambda}}_h,\mathbf{Ry}-\mathbf{y}_h)_{-1,1}\notag\\
		&\leq \norm{\mathbf{G}_h}_{-1} \norm{\nabla(\mathbf{Ry}-\mathbf{y}_h)}_{0,\Omega}-(\bm{\lambda}-\tilde{\bm{\lambda}}_h,\mathbf{Ry}-\mathbf{y}_h)_{-1,1}.\label{estimator_control}
	\end{align}
	Applying the Young's inequality in \eqref{estimator_control}, we get the following estimates:
	\begin{align}\label{estimator_control_1}
		\frac{\rho}{2}\norm{\nabla(\mathbf{Ry}-\mathbf{y}_h)}_{0,\Omega}^2 \leq &  \frac{1}{2\rho} \norm{\mathbf{G}_h}^2_{-1}-(\bm{\lambda}-\tilde{\bm{\lambda}}_h,\mathbf{Ry}-\mathbf{y}_h)_{-1,1}.
	\end{align}
	Also, using \eqref{galerkin_residual} one can easily derive
	\begin{align}\label{estimator_lambda}
		\norm{\bm{\lambda}-\tilde{\bm{\lambda}}_h}_{-1,\Omega}^2\leq (2+\frac{2}{\rho})\norm{\mathbf{G}_h}^2_{-1}-4(\bm{\lambda}-\tilde{\bm{\lambda}}_h,\mathbf{Ry}-\mathbf{y}_h)_{-1,1}.
	\end{align}
	Adding \eqref{estimator_control_1} and \eqref{estimator_lambda}, we get
\begin{align}\label{estimator_lambda_1}
	\norm{\nabla(\mathbf{Ry}-\mathbf{y}_h)}_{0,\Omega}^2+\norm{\bm{\lambda}-\tilde{\bm{\lambda}}_h}_{-1,\Omega}^2\leq (2+\frac{2}{\rho}+\frac{1}{\rho^2})\norm{\mathbf{G}_h}^2_{-1}-(4+\frac{2}{\rho})(\bm{\lambda}-\tilde{\bm{\lambda}}_h,\mathbf{Ry}-\mathbf{y}_h)_{-1,1}.
\end{align}
	From the equation \eqref{estimator_lambda_1} it is clear that the error in control and contact force densities are bounded by the dual norm of the Galerkin functional $\norm{\mathbf{G}_h}_{-1}$ and the duality pairing between the contact force densities and the controls.

\noindent 
First, we estimate the term $\norm{\mathbf{G}_h}_{-1}.$ Using the definition of $\mathbf{G}_h$ in \eqref{galerkin_residual}, with $\bm{\psi}\in \mathbf{Q}$, we obtain
	\begin{align}
		(\mathbf{G}_h,\bm{\psi})_{-1,1}=&\rho a(\mathbf{Ry}-\mathbf{y}_h,\bm{\psi})+(\bm{\lambda}-\tilde{\bm{\lambda}}_h,\bm{\psi})_{-1,1}\nonumber\\
		=&\rho a(\mathbf{Ry}-\mathbf{y}_h,\bm{\psi})+a(\bm{\psi},\bm{\phi}_h)-b(\bm{\psi},r_H)-(\mathbf{u}_h-\mathbf{u}_d,\bm{\psi})\nonumber\\&-\rho a(\mathbf{Ry},\bm{\psi})-(\tilde{\bm{\lambda}}_h,\bm{\psi})_{-1,1}\nonumber\\
		=&a(\bm{\psi},\bm{\phi}_h)-b(\bm{\psi},r_H)-(\mathbf{u}_h-\mathbf{u}_d,\bm{\psi})-\rho a(\mathbf{y}_h,\bm{\psi})-(\tilde{\bm{\lambda}}_h,\bm{\psi})_{-1,1}\nonumber\\
		=&(\mathbf{R}_h,\bm{\psi})_{-1,1}-(\tilde{\bm{\lambda}}_h,\bm{\psi})_{-1,1}\nonumber\\
		=& \sum_{i=1}^{2}(R_{h,i},\psi_i)_{-1,1}-\sum_{i=1}^{2}(\tilde{\lambda}_{h,i},\psi_i)_{-1,1}\nonumber\\
		=&\sum_{i=1}^{2}\sum_{p\in \cV_h}(R_{h,i},\psi_i\phi_p)_{-1,1}-\sum_{i=1}^{2}\sum_{p\in \cV_h^C}(\tilde{\lambda}^{p}_{h,i},\psi_i\phi_p)_{-1,1}\nonumber\\
		=&\sum_{i=1}^{2}\sum_{p\in \cV_h\setminus \cV_h^C}(R_{h,i},(\psi_i-m_p(\psi_i))\phi_p)_{-1,1}\nonumber\\&+\sum_{i=1}^{2}\sum_{p\in \cV_h^C}(R_{h,i},(\psi_i-m_p(\psi_i))\phi_p)_{-1,1}-\sum_{i=1}^{2}\sum_{p\in \cV_h^C}(\tilde{\lambda}^{p}_{h,i},\psi_i\phi_p)_{-1,1} 
		\nonumber\\
		=&\sum_{i=1}^{2}\sum_{p\in \cV_h\setminus \cV_h^C}(R_{h,i},(\psi_i-m_p(\psi_i))\phi_p)_{-1,1}\nonumber\\&+\sum_{i=1}^{2}\sum_{p\in \cV_h^C}(\tilde{R}_{h,i},(\psi_i-m_p(\psi_i))\phi_p)_{-1,1}\nonumber\\&-\sum_{i=1}^{2}\sum_{p\in \cV_h^C}\int_{\Gamma_{C}}\hat{\sigma}_i(\mathbf{y}_h)(\psi_i-m_p(\psi_i))\phi_p-\sum_{i=1}^{2}\sum_{p\in \cV_h^C}(\tilde{\lambda}^{p}_{h,i},\psi_i\phi_p)_{-1,1}
		\nonumber\\
		=&\sum_{i=1}^{2}\sum_{p\in \cV_h\setminus \cV_h^C}(R_{h,i},(\psi_i-m_p(\psi_i))\phi_p)_{-1,1}\nonumber\\&+\sum_{i=1}^{2}\sum_{p\in \cV_h^C}(\tilde{R}_{h,i},(\psi_i-m_p(\psi_i))\phi_p)_{-1,1}\nonumber\\&-\sum_{i=1}^{2}\sum_{p\in \cV_h^C\setminus \cV^{fC}_{h}}\int_{\Gamma_{C}}\hat{\sigma}_i(\mathbf{y}_h)(\psi_i-m_p(\psi_i))\phi_p.\label{lambda:exp}  		 
	\end{align}
	Here, we set $m_p(\psi_i)=0$ for Dirichlet and Neumann nodes. We exploited $(R_{h,i},\phi_p)_{-1,1}=0$ for all
	non-contact nodes and we inserted the definition
	of $\tilde{\bm{\lambda}}_h.$ Inserting definition of $R_{h,i}$ and $\tilde{R}_{h,i}$ in \eqref{lambda:exp}, we get  
	\begin{align}\label{exp:R_h}
		(\mathbf{G}_h,\bm{\psi})_{-1,1}=&\sum_{i=1}^{2}\sum_{p\in \cV_{h}} \int_{\omega_p}\mathcal{R}_i(\mathbf{y}_h)(\psi_i-m_p(\psi_i))\phi_p\nonumber\\&+\sum_{i=1}^{2}\sum_{p\in \cV_{h}}\int_{\gamma_{p,I}}\mathcal{J}^{I}_i(\mathbf{y}_h)(\psi_i-m_p(\psi_i))\phi_p\nonumber\\&-\sum_{i=1}^{2}\sum_{p\in \cV^C_{h}\setminus \cV^{fC}_{h}} \int_{\Gamma_C} \hat{\sigma}_i(\mathbf{y}_h)(\psi_i-m_p(\psi_i))\phi_p.
	\end{align}
where $\mathcal{R}_i(\mathbf{y}_h)=-\Delta \phi_{h,i}-\frac{\partial}{\partial x_i} r_H-(u_{h,i}-u_{d,i})+\rho \Delta y_{h,i}$ and 
$\mathcal{J}^{I}_i(\mathbf{y}_h)=\sjump{\nabla\phi_{h,i}+\begin{pmatrix} r_H \\ 0\end{pmatrix}'-\rho \nabla\mathbf{y}_h}.$
 For all nodes $p\in \cV_{h}\setminus(\cV^D_{h}\cup\cV^C_{h}\cup\cV^N_{h})$ we choose the constants
	\begin{align}\label{estm_c_p}
		m_p(\psi_i)=\frac{\int_{\omega_p}\psi_i\phi_p}{\int_{\omega_P}\phi_p}.
	\end{align}
	The mean value \eqref{estm_c_p} satisfies the following $L^2-$ approximation properties:
	\begin{align*}
		\norm{\psi_i-m_p(\psi_i)}_{0,\omega_p}&\leq C h_p \norm{\nabla \psi_i}_{0,\omega_p},\\
		\norm{\psi_i-m_p(\psi_i)}_{0,\gamma_p}&\leq C h_p^{\frac{1}{2}} \norm{\nabla \psi_i}_{0,\omega_p}\\
	\end{align*} 
	can be found in \cite[Lemma 3.1, Proposition 4.2]{AVesser_Verfurth2009}. For Dirichlet and Neumann nodes we have at least one edge $e\subset \gamma_{p,C}$, where the
	test function $\psi_i$ is zero, therefore we can deduce $\norm{\psi_i}_{0,\omega_p}\leq C h_p \norm{\nabla \psi_i}_{0,\omega_p}$ directly from
	the Poincar\'e-Friedrichs inequality. The above $L^2-$ approximation properties hold also for the
	constants $m_p(\psi_i)$ defined in \eqref{c_p_def_full_contact_a}, \eqref{c_p_def_full_contact_b} and \eqref{c_p_def_semi_contact} see \cite[Lemma 3.1, Proposition 4.2]{AVesser_Verfurth2009}.
	Using the above estimates and applying the Cauchy-Schwarz inequality in \eqref{exp:R_h} we arrive at
	\begin{align}\label{eq1234}
		(\mathbf{G}_h,\bm{\psi})_{-1,1}=&\sum_{i=1}^{2}\sum_{p\in \cV_{h}} h_p\norm{\mathcal{R}_i(\mathbf{y}_h)}_{0,\omega_p}\norm{\nabla\psi_i}_{0,\omega_p}\nonumber\\&+\sum_{i=1}^{2}\sum_{p\in \cV_{h}}h^{\frac{1}{2}}_p\norm{\mathcal{J}^{I}_i(\mathbf{y}_h)}_{0,\gamma_{p,I}}\norm{\nabla\psi_i}_{0,\omega_p}\nonumber\\&-\sum_{i=1}^{2}\sum_{p\in \cV^C_{h}\setminus \cV^{fC}_{h}}h_p^{\frac{1}{2}} \norm{\hat{\sigma}_i(\mathbf{y}_h)}_{0,\gamma_{p,C}} \norm{\nabla\psi_i}_{0,\omega_p}\nonumber\\&
		\lesssim \left(\sum_{p\in \cV_{h}} \big(h_p^2\norm{\mathcal{R}(\mathbf{y}_h)}^2_{0,\omega_p}+h_p \norm{\mathcal{J}^{I}(\mathbf{y}_h)}_{0,\gamma_{p,I}}^2\big)\right.\nonumber\\& \left.+\sum_{p\in \cV^C_{h}\setminus \cV^{fC}_{h}}h_p \norm{\hat{\bm{\sigma}}(\mathbf{y}_h)}_{0,\gamma_{p,C}}^2\right)^{\frac{1}{2}} \norm{\nabla\bm{\psi}}.
	\end{align}
	Thus from \eqref{eq1234}, we have 
	\begin{align}\label{R_h}
		\norm{\mathbf{G}_h}_{-1} \lesssim \left(\sum_{p\in \cV_{h}} \big(h_p^2\norm{\bm{\mathcal{R}}(\mathbf{y}_h)}^2_{0,\omega_p}+h_p \norm{\bm{\mathcal{J}}^{I}(\mathbf{y}_h)}_{0,\gamma_{p,I}}^2\big)+\sum_{p\in \cV^C_{h}\setminus \cV^{fC}_{h}}h_p \norm{\hat{\bm{\sigma}}(\mathbf{y}_h)}_{0,\gamma_{p,C}}^2\right)^{\frac{1}{2}}		
	\end{align}
where, $\bm{\mathcal{R}}(\mathbf{y}_h):=(\mathcal{R}_1(\mathbf{y}_h),\mathcal{R}_2(\mathbf{y}_h))=-\Delta \bm{\phi}_h-\nabla r_H-(\mathbf{u}_h-\mathbf{u}_d)+\rho \Delta \mathbf{y}_h$ and $\bm{\mathcal{J}}^{I}(\mathbf{y}_h):=(\mathcal{J}^{I}_1(\mathbf{y}_h),\mathcal{J}^{I}_2(\mathbf{y}_h))=\sjump{\nabla\bm{\phi}_h+r_H\I-\rho \nabla\mathbf{y}_h}.$

\noindent
Now, we need to find the upper bound of the term $(\tilde{\bm{\lambda}}_h-\bm{\lambda},\mathbf{Ry}-\mathbf{y}_h)_{-1,1}.$ We can write 
	\begin{align}\label{Estmlamda}
		(\tilde{\bm{\lambda}}_h-\bm{\lambda},\mathbf{Ry}-\mathbf{y}_h)_{-1,1}=(\tilde{\bm{\lambda}}_h,\mathbf{Ry}-\mathbf{y}_h)_{-1,1}+(\bm{\lambda},\mathbf{y}_h-\mathbf{Ry})_{-1,1}.
	\end{align}
	From the  equation \eqref{3.1e}, we have   $(\bm{\lambda},\mathbf{x}-\mathbf{Ry})_{-1,1}\leq 0$ for all $x\in\mathbf{Q}_{ad},$ taking $\mathbf{x}=\y_{h}$ we have $(\bm{\lambda},\mathbf{y}_h-\mathbf{Ry})_{-1,1}\leq 0.$ Thus it is enough to estimate the term  $(\tilde{\bm{\lambda}}_h,\mathbf{Ry}-\mathbf{y}_h)_{-1,1}.$ From the definition of the quasi discrete contact force density \eqref{Def_lambda_tilde}-\eqref{c_p_def_semi_contact}, we have
	\begin{align}\label{est:lambda_tilde}
		(\tilde{\bm{\lambda}}_h,\mathbf{Ry}-\mathbf{y}_h)_{-1,1}=&\sum_{i=1}^{2}\sum_{p\in \cV^C_{h}} (\tilde{\lambda}_{h,i}^p,(R_i\mathbf{y}-y_{h,i})\phi_p)_{-1,1}\nonumber\\
		=&\sum_{i=1}^{2}\left[\sum_{p\in \cV^{sC}_{h,a}}s_{p,i} m_p(R_i\mathbf{y}-y_{h,i})\int_{\gamma_{p,C}} \phi_p+\sum_{p\in \cV^{sC}_{h,b}}s_{p,i} m_p(R_i\mathbf{y}-y_{h,i})\int_{\gamma_{p,C}} \phi_p\right.\nonumber\\&\left.+\sum_{p\in \cV^{fC}_{h,a}}s_{p,i} m_p(R_i\mathbf{y}-y_{h,i})\int_{\gamma_{p,C}} \phi_p+\sum_{p\in \cV^{fC}_{h,b}}s_{p,i} m_p(R_i\mathbf{y}-y_{h,i})\int_{\gamma_{p,C}} \phi_p\right.\nonumber\\&\left.-\sum_{p\in \cV^{fC}_{h,a}}\int_{\gamma_{p,C}}\hat{\sigma}_i(\mathbf{y}_h) \big(R_i\mathbf{y}-y_{h,i}-m_p(R_i\mathbf{y}-y_{h,i})\big)\phi_p\right.\nonumber\\&\left.-\sum_{p\in \cV^{fC}_{h,b}}\int_{\gamma_{p,C}}\hat{\sigma}_i(\mathbf{y}_h) \big(R_i\mathbf{y}-y_{h,i}-m_p(R_i\mathbf{y}-y_{h,i})\big)\phi_p\right].
	\end{align} 
	We need to bound each term in the right hand side of \eqref{est:lambda_tilde}.\\ 
	\textit{First term:}
	\begin{align*}
		\sum_{p\in \cV^{sC}_{h,a}}s_{p,i} m_p(R_i\mathbf{y}-y_{h,i})\int_{\gamma_{p,C}} \phi_p&=\sum_{p\in \cV^{sC}_{h,a}}s_{p,i}\int_{\gamma_{p,C}} \phi_p \frac{\int_{\tilde{\gamma}_{p,C}}(R_i\mathbf{y}-y_{h,i})\phi_p}{\int_{\tilde{\gamma}_{p,C}}\phi_p}\\
		&=\sum_{p\in \cV^{sC}_{h,a}}s_{p,i}\int_{\gamma_{p,C}} \phi_p \frac{\int_{\tilde{\gamma}_{p,C}}(R_i\mathbf{y}-y^i_a)\phi_p}{\int_{\tilde{\gamma}_{p,C}}\phi_p}\\&+\sum_{p\in \cV^{sC}_{h,a}}s_{p,i}\int_{\gamma_{p,C}} \phi_p \frac{\int_{\tilde{\gamma}_{p,C}}(y^i_a-y_{h,i})\phi_p}{\int_{\tilde{\gamma}_{p,C}}\phi_p}\\
		&\leq \sum_{p\in \cV^{sC}_{h,a}}s_{p,i}\int_{\gamma_{p,C}} \phi_p \frac{\int_{\tilde{\gamma}_{p,C}}(y^i_a-y_{h,i})\phi_p}{\int_{\tilde{\gamma}_{p,C}}\phi_p}
		=\sum_{p\in \cV^{sC}_{h,a}} s_{p,i} d^{a}_{p,i},
	\end{align*}
	where $d^{a}_{p,i}=\int_{\tilde{\gamma}_{p,C}}(y^i_a-y_{h,i})\phi_p.$ In the above we exploit $s_{p,i}\leq 0,$  $R_i\mathbf{y}-y^i_a\geq 0$ and $\frac{\int_{\gamma_{p,C}}\phi_p}{\int_{\tilde{\gamma}_{p,C}}\phi_p}$ is a constant independent of $h_p$ if $\tilde{\gamma}_{p,C}$ is always a fixed fraction of ${\gamma}_{p,C}.$ 
	\par
	\noindent
	\textit{Second term:} 
	A similar arguments from the first term we have 
	\begin{align*}
		\sum_{p\in \cV^{sC}_{h,b}}s_{p,i} m_p(R_i\mathbf{y}-y_{h,i})\int_{\gamma_{p,C}} \phi_p&\leq \sum_{p\in \cV^{sC}_{h,b}} s_{p,i} d^{b}_{p,i},
	\end{align*}
	where $d^{b}_{p,i}=\int_{\tilde{\gamma}_{p,C}}(y^i_b-y_{h,i})\phi_p.$
	\par
	\noindent
	\textit{Third term:} For full contact nodes $p\in \cV^{fC}_{h,a},$ we have $R_i\mathbf{y}\geq y^i_a=y_{h,i}$ on $\gamma_{p,C}$ which implies $R_i\mathbf{y}-y_{h,i}\geq 0$ on $\gamma_{p,C}$ and therefore $m_p(R_i\mathbf{y}-y_{h,i})\geq 0.$ As further $s_{p,i}\leq 0$  we have $s_{p,i} m_p(R_i\mathbf{y}-y_{h,i})\leq 0.$
	Hence,
	\begin{align*}
		\sum_{p\in \cV^{fC}_{h,a}}s_{p,i} m_p(R_i\mathbf{y}-y_{h,i})\int_{\gamma_{p,C}} \phi_p\leq 0.
	\end{align*}
	
	\noindent
	\textit{Fourth term:}
	Similar to the third term for $p\in \cV^{fC}_{h,b}$ we have $s_{p,i}\geq0$
	and $R_i\mathbf{y}-y_{h,i}\leq0$ and hence
	$m_p(R_i\mathbf{y}-y_{h,i})\leq0.$
	Therefore 
	\begin{align*}
	\sum_{p\in \cV^{fC}_{h,b}}s_{p,i} m_p(R_i\mathbf{y}-y_{h,i})\int_{\gamma_{p,C}} \phi_p\leq 0.
	\end{align*}

	\par
	\noindent
	\textit{Fifth term:}
	\begin{align*}
		-&\int_{\gamma_{p,C}}\hat{\sigma}_i(\mathbf{y}_h) \big(R_i\mathbf{y}-y_{h,i}-m_p(R_i\mathbf{y}-y_{h,i})\big)\phi_p\\
		&=\sum_{e\subset \gamma_{p,C} }-\hat{\sigma}_i(\mathbf{y}_h)|_{e}\int_{e}\big(R_i\mathbf{y}-y_{h,i}-m_p(R_i\mathbf{y}-y_{h,i})\big)\phi_p\nonumber.
	\end{align*} 
	Here we have used the fact that $\hat{\sigma}_i(\mathbf{y}_h)|_{e}$ is piecewise constant on each edge. Since $p\in \cV^{fC}_{h,a}$ we have $-\hat{\sigma}_i(\mathbf{y}_h)|_{e}\leq 0$ and using the definition of $m_p$ from \eqref{c_p_def_full_contact_a} we get $$\int_{e}\big(R_i\mathbf{y}-y_{h,i}-m_p(R_i\mathbf{y}-y_{h,i})\big)\phi_p\geq 0.$$ Hence,
	\begin{align*}
		-\int_{\gamma_{p,C}}\hat{\sigma}_i(\mathbf{y}_h) \big(R_i\mathbf{y}-y_{h,i}-m_p(R_i\mathbf{y}-y_{h,i})\big)\phi_p\leq 0.
	\end{align*}  
	
	\par
	\noindent
	\textit{Sixth term:}
	Similar to the fourth term for $p\in \cV^{fC}_{h,b}$ we have $-\hat{\sigma}_i(\mathbf{y}_h)|_{e}\geq 0$ and using the definition of $m_p$ from \eqref{c_p_def_full_contact_b} we get $$\int_{e}\big(R_i\mathbf{y}-y_{h,i}-m_p(R_i\mathbf{y}-y_{h,i})\big)\phi_p\leq 0.$$ Hence,
	\begin{align*}
		-\int_{\gamma_{p,C}}\hat{\sigma}_i(\mathbf{y}_h) \big(R_i\mathbf{y}-y_{h,i}-m_p(R_i\mathbf{y}-y_{h,i})\big)\phi_p\leq 0.
	\end{align*}
	\noindent
	Using the all the above estimates from \textit{First term} to \textit{Sixth term}, in the right hand side of the equation \eqref{est:lambda_tilde}, we obtain
	\begin{align}\label{lambda}
		(\tilde{\bm{\lambda}}_h,\mathbf{Ry}-\mathbf{y}_h)_{-1,1}\leq & \sum_{i=1}^{2}\big[\sum_{p\in \cV^{sC}_{h,a}} s_{p,i} d^{a}_{p,i}+ \sum_{p\in \cV^{sC}_{h,b}} s_{p,i} d^{b}_{p,i}\big].
	\end{align}
	Substituting \eqref{R_h}, \eqref{Estmlamda}, and \eqref{lambda} in the right hand side of \eqref{estimator_lambda_1}, we obtain the following upper bound of control error:
	\begin{align}\label{est:control}
		\norm{\nabla(\mathbf{Ry}-\y_h)}_{0,\Omega}^2+&\norm{\bm{\lambda}-\tilde{\bm{\lambda}}_h}_{-1,\Omega}^2\lesssim 
		\sum_{p\in \cV_{h}} \big(\underbrace{h_p^2\norm{\bm{\mathcal{R}}(\mathbf{y}_h)}^2_{0,\omega_p}}_{\eta_\mathbf{y}^a}+\underbrace{h_p \norm{\bm{\mathcal{J}}^{I}(\mathbf{y}_h)}_{0,\gamma_{p,I}}^2}_{\eta_\mathbf{y}^b}\big)\nonumber\\&+\sum_{p\in \cV^C_{h}\setminus \cV^{fC}_{h}}\underbrace{h_p \norm{\hat{\bm{\sigma}}(\mathbf{y}_h)}_{0,\gamma_{p,C}}^2}_{\eta_\mathbf{y}^c}+\sum_{i=1}^{2}\big[\sum_{p\in \cV^{sC}_{h,a}}\underbrace{ s_{p,i} d^{a}_{p,i}}_{\eta_\mathbf{y}^d}+ \sum_{p\in \cV^{sC}_{h,b}} \underbrace{s_{p,i} d^{b}_{p,i}}_{\eta_\mathbf{y}^e}\big].
	\end{align}
	
	Denote the right hand side of \eqref{est:control} to be $\eta^2_\y.$ Thus substituting the estimates \eqref{apost_state}, \eqref{apost_adj_state}, and \eqref{est:control} in  \eqref{apostaux}, we prove Theorem \ref{aposteriori error estimator}.
\end{proof}
\begin{remark}
	The term $\eta_\mathbf{y}^e$ reminds of a complementarity condition. In fact, for a semi-contact node $s_{p,i}d_{p,i}^b$ would be a complementarity condition with respect to the quasi-discrete contact force density $(\tilde{\lambda}_{h,i}^p,(y_b^i-y_{h,i})\phi_p)_{-1,1}$ if $\tilde{\gamma}_{p,C}$ was replaced by $\gamma_{p,C}$. Thus we refer to $\eta_\mathbf{y}^d$ and $\eta_\mathbf{y}^e$ as complementarity residual and call $\eta_\mathbf{y}^c$ contact
	stress residual. The contributions $\eta_\mathbf{y}^c,$ $\eta_\mathbf{y}^d$ are localized to semi-contact nodes and nodes which are not actually in contact. In the unconstrained case, we have $\eta_\mathbf{y}^d=\eta_\mathbf{y}^e=0$
	and $\eta_\mathbf{y}^c$ has contributions from all potential contact nodes such that $\eta_\mathbf{y}$ is a residual
	error estimator for linear elliptic boundary value problems where the potential contact	boundary is replaced by a Neumann boundary with zero Neumann data.
\end{remark}

\begin{theorem}[\bf Local Efficiency]\label{local_efficiency}
	Let $\mathcal{T}_{e}$ be the set of two triangles sharing the edge $e\in \mathcal{E}_{h}^{i}.$ Then, it hold
	\begin{align}
		h_T\lVert \mathbf{f}+\Delta \mathbf{u}_h-\nabla p_H\rVert_{0,T} \lesssim &\big(\lVert \nabla(\mathbf{u}-\mathbf{u}_{h})\rVert_{0,T}+\lVert p-p_{H}\rVert_{0,T}\nonumber\\&+ osc(\mathbf{f},T)\big),\nonumber\\
		h_{T} \lVert \Delta \bm{\phi}_h+\nabla r_H+ \mathbf{u}_{h}-\mathbf{u}_d\rVert_{0,T}\lesssim & \big(\lVert \nabla(\mathbf{u}-\mathbf{u}_{h})\rVert_{0,T}+\lVert \nabla( \bm{\phi}-\bm{\phi}_{h})\rVert_{0,T}\nonumber\\&+\lVert r-r_{H}\rVert_{0,T}+ osc(\mathbf{u}_d,T)\big),\nonumber\\
		h_p\norm{\bm{\mathcal{R}}(\mathbf{y}_h)}_{0,\omega_p} \lesssim & \big(\lVert \nabla(\mathbf{u}-\mathbf{u}_{h})\rVert_{0,\omega_p}+\rho\lVert \nabla(\mathbf{y}-\mathbf{y}_{h})\rVert_{0,\omega_p}\nonumber\\&+\lVert \nabla( \bm{\phi}-\bm{\phi}_{h})\rVert_{0,\omega_p}+\lVert r-r_{H}\rVert_{0,\omega_p}\nonumber\\&+ osc(\mathbf{u}_d,\omega_p)\big)\label{Eff:Rh}\\
		\lVert h_e^{\frac{1}{2}}\sjump{p_{H}I-\nabla{\mathbf{u}_{h}}}\rVert_{0,e}\lesssim &  \sum_{T\in \mathcal{T}_{e}}\big(\lVert \nabla(\mathbf{u}-\mathbf{u}_{h})\rVert_{0,T}+\lVert p-p_{H}\rVert_{0,T}\nonumber\\&+ osc(\mathbf{f},T)\big),\nonumber\\
		\lVert h^{\frac{1}{2}}\sjump{r_{H}I+\nabla{\bm{\phi}_{h}}}\rVert_{0,e}\lesssim &  \sum_{T\in\mathcal{T}_{e}}\big(\lVert \nabla(\mathbf{u}-\mathbf{u}_{h})\rVert_{0,T}+\lVert \nabla( \bm{\phi}-\bm{\phi}_{h})\rVert_{0,T}\nonumber\\&+\lVert r-r_{H}\rVert_{0,T}+ osc(\mathbf{u}_d,T)\big),\nonumber\\
		\lVert h_p^{\frac{1}{2}}\bm{\mathcal{J}}^{I}(\mathbf{y}_h)\rVert_{0,\gamma_{p,I}}\lesssim &  \lVert \nabla(\mathbf{u}-\mathbf{u}_{h})\rVert_{0,T}+\lVert \nabla( \bm{\phi}-\bm{\phi}_{h})\rVert_{0,\omega_p}\nonumber\\&+\lVert r-r_{H}\rVert_{0,\omega_p}+\rho\lVert \nabla(\mathbf{y}-\mathbf{y}_{h})\rVert_{0,\omega_p}\nonumber\\&+ osc(\mathbf{u}_d,\omega_p),\label{Eff:JI}\\
		\norm{\nabla \cdot \mathbf{u}_h}_{0,T}\lesssim &  \norm{\nabla(\mathbf{u}-\mathbf{u}_h)}_{0,T}\nonumber\\
		\norm{\nabla \cdot \bm{\phi}_h}_{0,T}\lesssim &  \norm{\nabla(\bm{\phi}-\bm{\phi}_h)}_{0,T}.\nonumber
	\end{align}
	Further, for any Neumann boundary edge $e\in \mathcal{E}_{h}^{b,N}$, it hold
	\begin{align*}
		\lVert h_e^{\frac{1}{2}}\sjump{p_{H}I-\nabla{\mathbf{u}_{h}}}\rVert_{0,e}\lesssim&  \big(\lVert \nabla(\mathbf{u}-\mathbf{u}_{h})\rVert_{0,T}+\lVert p-p_{H}\rVert_{0,T}+ osc(\mathbf{f},T)\big),\\
		\lVert h_e^{\frac{1}{2}}\sjump{r_{H}I+\nabla{\bm{\phi}_{h}}}\rVert_{0,e}\lesssim &  \big(\lVert \nabla ( \mathbf{u}-\mathbf{u}_{h})\rVert_{0,T}+\lVert \nabla( \bm{\phi}-\bm{\phi}_{h})\rVert_{0,T}+\lVert r-r_{H}\rVert_{0,T}\\&+ osc(\mathbf{u}_d,T)\big).
	\end{align*}
\end{theorem}
In the contact zone, for $p\in \cV^C_{h}\setminus \cV^{fC}_{h}$
\begin{align}\label{efficiency_sigma1}
	 h^{\frac{1}{2}}_p \lVert\hat{\bm{\sigma}}(\mathbf{y}_h)\rVert&_{0,\gamma_{p,C}} \lesssim \big( \lVert \nabla (\mathbf{u}-\mathbf{u}_{h})\rVert_{0,\omega_p}+\lVert \nabla( \bm{\phi}-\bm{\phi}_{h})\rVert_{0,\omega_p}\nonumber\\&+\lVert r-r_{H}\rVert_{0,\omega_p}+\rho\lVert \nabla(\mathbf{y}-\mathbf{y}_{h})\rVert_{0,\omega_p}+\norm{\bm{\lambda}-\tilde{\bm{\lambda}}_h}_{-1,\omega_p}\nonumber\\&+ osc(\mathbf{u}_d,\omega_p)\big),  
\end{align}
for $p\in \cV^{sC}_{h,a}$
\begin{align}\label{efficiency_s_p^1}
	 \sum_{i=1}^{2}s_{p,i} d^{a}_{p,i} \lesssim &  \big( \lVert \nabla (\mathbf{u}-\mathbf{u}_{h})\rVert_{0,\omega_p}+\lVert \nabla( \bm{\phi}-\bm{\phi}_{h})\rVert_{0,\omega_p}+\lVert r-r_{H}\rVert_{0,\omega_p}\nonumber\\&+\rho\lVert \nabla(\mathbf{y}-\mathbf{y}_{h})\rVert_{0,\omega_p}+\norm{\bm{\lambda}-\tilde{\bm{\lambda}}_h}_{-1,\omega_p}+ osc(\mathbf{u}_d,\omega_p)\big)^2
\end{align}
for $p\in \cV^{sC}_{h,b}$
\begin{align}\label{efficiency_s_p}
	\sum_{i=1}^{2} s_{p,i} d^{b}_{p,i} \lesssim &  \big( \lVert \nabla (\mathbf{u}-\mathbf{u}_{h})\rVert_{0,\omega_p}+\lVert \nabla( \bm{\phi}-\bm{\phi}_{h})\rVert_{0,\omega_p}+\lVert r-r_{H}\rVert_{0,\omega_p}\nonumber\\&+\rho\lVert \nabla(\mathbf{y}-\mathbf{y}_{h})\rVert_{0,\omega_p}+\norm{\bm{\lambda}-\tilde{\bm{\lambda}}_h}_{-1,\omega_p}+ osc(\mathbf{u}_d,\omega_p)\big)^2
\end{align}
where oscillation of a given function $\mathbf{f}, \mathbf{u}_d\in \mathbf{L}^2(T)$ is defined by $$osc(\mathbf{f},T)=h_T \min_{\mathbf{f}_h\in\mathbf{P}_0(T)} \norm{\mathbf{f}-\mathbf{f}_h}_{0,T}, osc(\mathbf{u}_d,T)=h_T \min_{\mathbf{g}_h\in\mathbf{P}_0(T)} \norm{\mathbf{u}_d-\mathbf{g}_h}_{0,T},$$ 
and similarly, we define the oscillation, $osc(\mathbf{u}_d,\omega_p)=h_T \min_{\mathbf{g}_h\in\mathbf{P}_0(\omega_p)} \norm{\mathbf{u}_d-\mathbf{g}_h}_{0,\omega_p}.$
\begin{proof}
	The local efficiencies in the above theorem can be deduced by the standard bubble function techniques in \cite{Verfurth:1995}, except the terms \eqref{efficiency_sigma1}, \eqref{efficiency_s_p^1}, and \eqref{efficiency_s_p}. First, we will prove the efficiency \eqref{efficiency_sigma1}. We make use of the relation
	between the Galerkin functional and the quantity of interest which here is the boundary stress. It directly follows from the definition of the Galerkin functional \eqref{galerkin_residual}, \eqref{eq:VI} and \eqref{3.1e} that
	\begin{align}
		\norm{\mathbf{G}_h}_{-1,\omega_p}&\lesssim \norm{\nabla(\mathbf{Ry}-\mathbf{y}_h)}_{0,\omega_p}+\norm{\bm{\lambda}-\tilde{\bm{\lambda}}_h}_{-1,\omega_p}\nonumber\\
		&\lesssim \lVert \nabla(\mathbf{y}-\mathbf{y}_{h})\rVert_{0,\omega_p}+\lVert \nabla (\mathbf{u}-\mathbf{u}_{h})\rVert_{0,\omega_p}+\lVert \nabla( \bm{\phi}-\bm{\phi}_{h})\rVert_{0,\omega_p}\nonumber\\&+\lVert r-r_{H}\rVert_{0,\omega_p}+\norm{\bm{\lambda}-\tilde{\bm{\lambda}}_h}_{-1,\omega_p}.\label{Sto:DC:estmGh}
	\end{align}
	Let $\bar{p}\in \cV^C_{h}\setminus \cV^{fC}_{h}$ be an arbitrary but fixed node. In the following $s$ denotes a side
	which belongs to $\gamma_{\bar{p},C}$. We take the corresponding side bubble function
	$\xi_s:=\prod_{p\in s}\phi_p.$
	Test the function $\xi_s\mathbf{e}_1$ in the equation \eqref{lambda:exp}, we get
	\begin{align}
		\sum_{p\in \cV_h^C\setminus \cV^{fC}_{h}}\int_{\Gamma_{C}}\hat{\sigma}_1(\mathbf{y}_h)\xi_s\phi_p=&-(\mathbf{G}_h,\xi_s\mathbf{e}_1)_{-1,1}+\sum_{p\in \cV_h\setminus \cV_h^C}(R_{h,1},\xi_s\phi_p)_{-1,1}\nonumber\\&+\sum_{p\in \cV_h^C}(\tilde{R}_{h,1},(\xi_s-m_p(\xi_s))\phi_p)_{-1,1}+\sum_{p\in \cV_h^C\setminus \cV^{fC}_{h}}\int_{\Gamma_{C}}\hat{\sigma}_1(\mathbf{y}_h)m_p(\xi_s)\phi_p\nonumber\\&
		=-(\mathbf{G}_h,\xi_s\mathbf{e}_1)_{-1,1}+\sum_{p\in \cV_{h}} \int_{\omega_s}\mathcal{R}_1(\mathbf{y}_h)\xi_s\phi_p\nonumber\\&-\sum_{p\in \cV_h^{sC}}s_{p,1} m_p(\xi_s)\int_{\gamma_{p,C}}\phi_p-\sum_{p\in \cV_h^{fC}}(\tilde{R}_{h,1},\phi_p)_{-1,1}m_p(\xi_s)\label{sto_dc_500}	
	\end{align}
	If the side $s$ is not contained in any patch $\gamma_{p,C}$ of semi- or full-contact nodes $p$, the
	two last terms are zero and we can proceed similar to the case of \eqref{Eff:Rh} and \eqref{Eff:JI}.
	Otherwise, in order to get rid of the last two terms, we replace $\xi_s$ by a suitable function
	$\theta_s$ such that $m_p(\theta_s)=0$ for all semi- and full-contact nodes. The value of $m_p(.)$
	for a semi-contact node $p$ depends on $\tilde{\gamma}_{p,C}$ which is a strict subset of $\gamma_{p,C}$ compare \eqref{c_p_def_semi_contact}. If $\gamma_{p,C}$ consists of two intervals we choose the inner third of $\gamma_{p,C}$ containing $p$ as $\tilde{\gamma}_{p,C}.$
	
	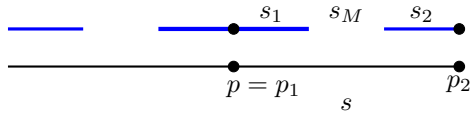
\begin{figure}[h!]
		\begin{center}
			\begin{tikzpicture}
				\draw[black,thick] (-3,0) -- (3,0);
				\draw[blue,very thick] (-3,.5) -- (-2,.5);
				\draw[blue,ultra thick] (-1,.5) -- (1,.5);
				\draw[blue,very thick] (2,.5) -- (3,.5);
				\filldraw[black] (0,0) circle (2pt); 
				\filldraw[black] (3,0) circle (2pt) node[anchor=north] {$p_2$};
				\filldraw[black] (3,.5) circle (2pt);
				\filldraw[black] (0,.5) circle (2pt);
				\node[black,thick] () at (1.5,-.5){$s$};
				\node[black,thick] () at (.4,-.3){$p=p_1$};
				\node[black,thick] () at (.5,.7){$s_1$};
				\node[black,thick] () at (2.5,.7){$s_2$};
				\node[black,thick] () at (1.5,.7){$s_M$};
			\end{tikzpicture}
			\caption{ Subgrid of boundary patch $\gamma_{p,C}$}
			\label{Subgrid_of_boundary_patch}
		\end{center}
	\end{figure}
For example
in Fig. \ref{Subgrid_of_boundary_patch}, the dark blue region is $\tilde{\gamma}_{p,C}$ for $p=p_1$. A side $s$ has two nodes $p_1,p_2$.
We denote the sides of the subgrid containing $p_i$ by $s_i$ and the middle part by $s_M$, see
Fig. \ref{Subgrid_of_boundary_patch}.
For the function $\theta_s$ we make the following ansatz
\begin{align}\label{ansatz}
	\theta_s=\sum_{i=1}^{2}a_i\xi_i+a_M\xi_M,
\end{align}	
where $\xi_i$ and $\xi_M$ are side bubble functions to $s_i$ and $s_M$. The coefficients $a_1,a_2,a_M$ are determined so that $\int_{s}1=\sum_{p\in \cV_h^C\setminus \cV^{fC}_{h}}\int_{s}\theta_s\phi_p,$ and $\int_{s_i}\theta_s\phi_{p_i}=0$ for $p_i$ full-contact or semi-contact nodes. As $\bar{p}$ is not a full-contact node there is at least one contribution in the right hand side of the first condition. Inserting the ansatz \eqref{ansatz} in the aforementioned conditions,
we get a solvable system of three equations with three coefficients (degrees of freedom) $a_1,a_2,a_M$. At this point the special choice of $	m_p(\phi)=\frac{\int_{\tilde{\gamma}_{p,C}}\phi\phi_p}{\int_{\tilde{\gamma}_{p,C}}\phi_p}$ 
as mean value on $\tilde{\gamma}_{p,C}$ for semi-contact nodes becomes important because
the choice
$	m_p(\phi)=\frac{\int_{\gamma_{p,C}}\phi\phi_p}{\int_{\gamma_{p,C}}\phi_p}$
as mean value over the whole patch $\gamma_{p,C}$ would lead to a
contradiction of the conditions. In the second condition $s_i$ would be replaced by $s$ and
the condition $\int_{s}\theta_s\phi_{p_i}=0$ for all $p_i$ of the side $s$ would imply
$\sum_{p\in \cV_h^C\setminus \cV^{fC}_{h}}\int_{s}\theta_s\phi_p=0$
such that the first condition could not be fulfilled.
As we assumed that the mesh is made of simplices, $\hat{\sigma}_1(\mathbf{y}_h)$ is constant on $s$. Consequently, $m_p(\theta_s)=0$ implies $m_p(\hat{\sigma}_1(\mathbf{y}_h)\theta_s)=0$ and it follows from the first condition	
\begin{align}\label{sigma_hat}
	\norm{\hat{\sigma}_1(\mathbf{y}_h)}_{0,s}^2=\sum_{p\in \cV_h^C\setminus \cV^{fC}_{h}}\int_{s}\hat{\sigma}_1(\mathbf{y}_h)\hat{\sigma}_1(\mathbf{y}_h)\theta_s\phi_p.
\end{align}	
Putting together \eqref{sigma_hat}, \eqref{sto_dc_500} with test function $\hat{\sigma}_1(\mathbf{y}_h)\theta_s$  instead of $\xi_s$ and exploiting
the conditions $m_p(\hat{\sigma}_1(\mathbf{y}_h)\theta_s)=0$ for all contact nodes we end up with
\begin{align}
	\norm{\hat{\sigma}_1(\mathbf{y}_h)}_{0,s}^2&=\sum_{p\in \cV_h^C\setminus \cV^{fC}_{h}}\int_{s}\hat{\sigma}_1(\mathbf{y}_h)\hat{\sigma}_1(\mathbf{y}_h)\theta_s\phi_p\nonumber\\
	&= -(\mathbf{G}_h,\hat{\sigma}_1(\mathbf{y}_h)\theta_s\mathbf{e}_1)_{-1,1}+ \int_{\omega_s}\mathcal{R}_1(\mathbf{y}_h)\hat{\sigma}_1(\mathbf{y}_h)\theta_s\nonumber\\&
	\lesssim \norm{\mathbf{G}_h}_{-1,\omega_{\bar{p}}}\norm{\hat{\sigma}_1(\mathbf{y}_h)\theta_s}_{1,\omega_s}+\norm{\mathcal{R}_1(\mathbf{y}_h)}_{0,\omega_s}\norm{\hat{\sigma}_1(\mathbf{y}_h)\theta_s}_{1,\omega_s}\nonumber\\&
	\lesssim \norm{\mathbf{G}_h}_{-1,\omega_{\bar{p}}}h_s^{-\frac{1}{2}}\norm{\hat{\sigma}_1(\mathbf{y}_h)}_{0,s}+h_s^{\frac{1}{2}}\norm{\mathcal{R}_1(\mathbf{y}_h)}_{0,\omega_s}\norm{\hat{\sigma}_1(\mathbf{y}_h)}_{0,s}\label{stoDC:550}
\end{align}	
where $h_s:= diam(s).$ In the last line of \eqref{stoDC:550} we used the properties of the bubble
functions on the subgrid and the fact that $\tilde{\gamma}_{p,C}$ is a fixed portion of $\gamma_{p,C}$ so that $h_s=c h_{s_i}$
for a mesh-independent constant $c$. We divide by $h_s^{-\frac{1}{2}}\norm{\hat{\sigma}_1(\mathbf{y}_h)}_{0,s}$ leading
to 
\begin{align}
	h_s^{\frac{1}{2}}\norm{\hat{\sigma}_1(\mathbf{y}_h)}_{0,s}\lesssim \norm{\mathbf{G}_h}_{-1,\omega_{\bar{p}}}+h_s\norm{\bm{\mathcal{R}}(\mathbf{y}_h)}_{0,\omega_s}.
\end{align}
By means of the triangle inequality, the shape-regularity, $h_s\approx h_p$ and the upper bounds \eqref{Sto:DC:estmGh} 
and \eqref{Eff:Rh} of $\norm{\mathbf{G}_h}_{-1,\omega_{\bar{p}}}$and $\norm{\bm{\mathcal{R}}(\mathbf{y}_h)}_{0,\omega_s},$ we get
\begin{align}\label{sigma_1_estm}
	h_s^{\frac{1}{2}}\norm{\hat{\sigma}_1(\mathbf{y}_h)}_{0,\gamma_{p,C}}&\lesssim \lVert \nabla (\mathbf{u}-\mathbf{u}_{h})\rVert_{0,\omega_p}+\lVert \nabla( \bm{\phi}-\bm{\phi}_{h})\rVert_{0,\omega_p}+\lVert r-r_{H}\rVert_{0,\omega_p}\nonumber\\&+\rho\lVert \nabla(\mathbf{y}-\mathbf{y}_{h})\rVert_{0,\omega_p}+\norm{\bm{\lambda}-\tilde{\bm{\lambda}}_h}_{-1,\omega_p}+ osc(\mathbf{u}_d,\omega_p).
\end{align}
Similarly, one can derive 
\begin{align}\label{sigma_2_estm}
	h_s^{\frac{1}{2}}\norm{\hat{\sigma}_2(\mathbf{y}_h)}_{0,\gamma_{p,C}}&\lesssim \lVert \nabla (\mathbf{u}-\mathbf{u}_{h})\rVert_{0,\omega_p}+\lVert \nabla( \bm{\phi}-\bm{\phi}_{h})\rVert_{0,\omega_p}+\lVert r-r_{H}\rVert_{0,\omega_p}\nonumber\\&+\rho\lVert \nabla(\mathbf{y}-\mathbf{y}_{h})\rVert_{0,\omega_p}+\norm{\bm{\lambda}-\tilde{\bm{\lambda}}_h}_{-1,\omega_p}+ osc(\mathbf{u}_d,\omega_p)
\end{align}

Adding \eqref{sigma_1_estm} and \eqref{sigma_2_estm}, we get the desired result \eqref{efficiency_sigma1}.

	Now we turn back to the terms \eqref{efficiency_s_p} and \eqref{efficiency_s_p^1}.
	   The proof of both \eqref{efficiency_s_p} and \eqref{efficiency_s_p^1} proceeds similarly so we will give a sketch of the proof of \eqref{efficiency_s_p} here and details can be found in\cite[Sec. 5.2]{Andreas_Veeser_Signorini:2015}. To prove \eqref{efficiency_s_p}, we derive a lower bound of the local error in terms of the local contributions of $s_{p,1} d^{b}_{p,1}.$ If $s_{p,1}=0$ or $(y^1_b-y_{h,1})(q)= 0$ for all neighbouring nodes of $p$ we
	have $s_{p,1}d^{b}_{p,1}=0$. Therefore, we assume $s_{p,1}>0$ and $(y^1_b-y_{h,1})(q)>0$ for at least one
	node on $\gamma_{p,C}$. Let $\hat{q}$ be a node which fulfills $(y^1_b-y_{h,1})(\hat{q})\geq(y^1_b-y_{h,1})(q)$ for all
	neighboring nodes $q$ of $p$. Due to $s_{p,1}>0$ we have $(y^1_b-y_{h,1})(p) = 0.$ As we consider
	boundary meshes of triangles and intervals the discrete functions are piecewise linear.
	Using Taylor series expansion of $(y^1_b-y_{h,1})$ about $p,$ we get 
	\begin{align}\label{y_b-y_h}
		(y^1_b-y_{h,1})(\hat{q})= \nabla|_{\hat{e}}(y^1_b-y_{h,1})\cdot (\hat{q}-p)\lesssim h_p \nabla|_{\hat{e}}(y^1_b-y_{h,1})\cdot \tau
	\end{align}
	where $\hat{e}\subset \gamma_{p,C}$ is an edge containing the nodes $\hat{q}$ and $p$ and $\tau$ is the unit tangential vector pointing from $p$ to $\hat{q}.$ The following estimate of \eqref{y_b-y_h} follows from \cite[Sec. 5.2]{Andreas_Veeser_Signorini:2015}:
	\begin{align*}
		(y^1_b-y_{h,1})(\hat{q})
		& \lesssim h_p^{\frac{1}{2}}\lVert \bm{\mathcal{J}}^{I}(\mathbf{y}_h)\rVert_{0,\gamma_{p,I}}.
	\end{align*}
	Since, $(y^1_b-y_{h,1})(q)\leq (y^1_b-y_{h,1})(\hat{q})$ for all $q\in\gamma_{p,C}$ , we can conclude that 
	\begin{align}\label{d^{b}_p expression}
		d^{b}_{p,1}=\int_{\tilde{\gamma}_{p,C}}(y^1_b-y_{h,1})\phi_p\lesssim h^{\frac{3}{2}}_p \lVert \bm{\mathcal{J}}^{I}(\mathbf{y}_h)\rVert_{0,\gamma_{p,I}}.
	\end{align}
	Now,
	\begin{align*}
		s_{p,1}d^{b}_{p,1}&= \frac{\int_{\omega_p}\mathcal{R}_1(\mathbf{y}_h) \phi_p+\int_{\gamma_{p,I}}  \mathcal{J}^{I}_1(\mathbf{y}_h)\phi_p-\int_{\gamma_{p,C}}\hat{\sigma}_1(\mathbf{y}_h)\phi_p}{\int_{\gamma_{p,C}}\phi_p} \cdot \int_{\tilde{\gamma}_{p,C}}(y^1_b-y_{h,1})\phi_p\\&
		\lesssim \big(h_p\norm{\mathcal{R}_1(\mathbf{y}_h)}_{0,\omega_p} +h_p^{\frac{1}{2}}  \norm{\mathcal{J}^{I}_1(\mathbf{y}_h)}_{0,\gamma_{p,I}}+h_p^{\frac{1}{2}}\norm{\hat{\sigma}_1(\mathbf{y}_h)}_{0,\gamma_{p,C}}\big) \cdot h^{-1}_p\cdot \int_{\tilde{\gamma}_{p,C}}(y^1_b-y_{h,1})\phi_p\\&
		\lesssim \big(h_p\norm{\bm{\mathcal{R}}(\mathbf{y}_h)}_{0,\omega_p} +h_p^{\frac{1}{2}}  \norm{\bm{\mathcal{J}}^{I}(\mathbf{y}_h)}_{0,\gamma_{p,I}}+h_p^{\frac{1}{2}}\norm{\hat{\bm{\sigma}}(\mathbf{y}_h)}_{0,\gamma_{p,C}}\big)\cdot h^{\frac{1}{2}}_p \lVert \bm{\mathcal{J}}^{I}(\mathbf{y}_h)\rVert_{0,\gamma_{p,I}}
	\end{align*}
	Applying Cauchy-Schwarz inequality, \eqref{efficiency_sigma1}, \eqref{Eff:Rh} and \eqref{Eff:JI} we obtain
	\begin{align*}
		s_{p,1}d^{b}_{p,1}\lesssim& \big(\lVert \nabla (\mathbf{u}-\mathbf{u}_{h})\rVert_{0,\omega_p}+\lVert \nabla( \bm{\phi}-\bm{\phi}_{h})\rVert_{0,\omega_p}+\lVert r-r_{H}\rVert_{0,\omega_p}\\&+\rho\lVert \nabla(\mathbf{y}-\mathbf{y}_{h})\rVert_{0,\omega_p}+\norm{\bm{\lambda}-\tilde{\bm{\lambda}}_h}_{-1,\omega_p}+ osc(\mathbf{u}_d,\omega_p)\big)^2.
	\end{align*}
	Thus, one can obtain the upper bound of $s_{p,2}d^{b}_{p,2}.$ Hence, for $p\in \cV^{sC}_{h,b}$ we have
	\begin{align*}
		\sum_{i=1}^{2}s_{p,i} d^{b}_{p,i} \lesssim &  \big( \lVert \nabla (\mathbf{u}-\mathbf{u}_{h})\rVert_{0,\omega_p}+\lVert \nabla( \bm{\phi}-\bm{\phi}_{h})\rVert_{0,\omega_p}+\lVert r-r_{H}\rVert_{0,\omega_p}\nonumber\\&+\rho\lVert \nabla(\mathbf{y}-\mathbf{y}_{h})\rVert_{0,\omega_p}+\norm{\bm{\lambda}-\tilde{\bm{\lambda}}_h}_{-1,\omega_p}+ osc(\mathbf{u}_d,\omega_p)\big)^2.
	\end{align*}	
Similarly, we can prove \eqref{efficiency_s_p^1}.
\end{proof}

\section{Numerical Experiments}\label{sec:Numerics}
The aim of the given section is to numerically illustrate the theoretical results derived in Sections \ref{Discrete Problems} and \ref{Sto:DC_error_analysis}, respectively. We conduct two experiments with two model problems, one is a smooth solution on square mesh the other is a non-smooth solution on a non-convex domain. We construct the model problems with known solutions. The numerical experiments are performed on two model problems using MATLAB(version R2021a) software. For the computational simplicity, we slightly modify the cost functional $J$, denoted by $\tilde{J}$, by

\begin{align}\label{Modified_prob}
\tilde{J}(\mathbf{w},\mathbf{x})=\frac{1}{2}\|\mathbf{w}-\mathbf{u}_d\|^2+\frac{\rho}{2}\|\nabla (\mathbf{x}-\mathbf{y}_d)\|^2 
\end{align}

subject to PDE,
\begin{equation}\label{modified_eq}
	\begin{split}
	-\Delta \mathbf{w}+\nabla{p}&=\mathbf{f} \quad \text{in}\;\Omega,\\
	\nabla\cdot{\mathbf{w}}&=0 \quad \text{in}\; \Omega,\\
	\bf{w}&=\mathbf{x} \quad \text{on}\; \Gamma_C,\\
	\bf{w}&=\mathbf{0} \quad \text{on}\; \Gamma_D, 
	\end{split}
\end{equation}
the set of controls is given by,
$$\mathbf{Q}_{ad}:=\{\mathbf{x}\in \mathbf{H}_D^1(\Omega): \mathbf{y}_a\leq\bm{\gamma}_{0}(\mathbf{x})\leq \mathbf{y}_b\text{ a.e. on } \Gamma_C\},$$

\noindent
where the space $\mathbf{H}^1_D(\Omega)$ consists of $\mathbf{H}^1(\Omega)$ functions with vanishing trace on $\Gamma_D.$ The function $\mathbf{y}_d$ is given, and $\partial\Omega=\Gamma_C\cup \bar{\Gamma}_D$. Then the minimization problem reads: Find
 $(\mathbf{u},\mathbf{y})\in \mathbf{H}_{D}^1(\Omega) \times \mathbf{Q}_{ad}$ satisfies \eqref{modified_eq} such that
 \begin{align*}
 \tilde{J}(\mathbf{u},\mathbf{y})=\min_{(\mathbf{v},\mathbf{x})\in \mathbf{H}_{D}^1(\Omega)\times \mathbf{Q}_{ad}} \tilde{J}(\mathbf{v},\mathbf{x}).
 \end{align*}

 The corresponding discrete optimality system is given by
 \begin{equation}\label{modified_discModel}
 	\begin{split}
 		\mathbf{u}_h&=\mathbf{w}_h+\mathbf{y}_h,\quad \mathbf{w}_h\in \mathbf{V}_h,\\
 		a(\mathbf{w}_h,\mathbf{z}_h)+b(\mathbf{z}_h,p_H) &=( \mathbf{f},\mathbf{z}_h)-a(\mathbf{y}_h,\mathbf{z}_h) \;\;\;
 		{\rm for~all}\;\mathbf{z}_h \in \mathbf{V}_h,\\ 
 		b(\mathbf{u}_h,q_H)&=0\;\quad {\rm for~all}\; q_H \in M_H ,\\
 		a(\mathbf{z}_h,\bm{\phi}_h)-b(\mathbf{z}_h,r_H) &=( \mathbf{u}_h-\mathbf{u_d},\mathbf{z}_h) \;\;\;
 		{\rm for~all}\;\mathbf{z}_h \in \mathbf{V}_h ,\\ 
 		b(\bm{\phi}_h,q_H)&=0\;\quad {\rm for~all}\; q_H \in M_H ,\\
 		\rho\, a(\mathbf{y}_h,\mathbf{x}_h-\mathbf{y}_h)\geq& a(\mathbf{x}_h-\mathbf{y}_h,\bm{\phi}_h)-b(\mathbf{x}_h-\mathbf{y}_h,r_H)\\&-(\mathbf{u}_h-\mathbf{u_d},\mathbf{x}_h-\mathbf{y}_h)+\rho a(\mathbf{y}_d,\mathbf{x}_h-\mathbf{y}_h)\;\quad {\rm for~all}\; \mathbf{x}_h\in
 		\mathbf{Q}^h_{ad},
 	\end{split}
 \end{equation}where, 
 $\mathbf{V}_h:=\{\mathbf{v}_h\in \mathbf{H}_{0}^1(\Omega): \mathbf{v}_h|_T\in \mathbf{P}_1(T)\;\; \forall T \in \mathcal{T}_h\}$
 and $\mathbf{Q}^h_{ad}=\mathbf{Q}_h\cap\mathbf{Q}_{ad}.$ The set $\mathbf{Q}_h$ is defined by
 $\mathbf{Q}_h=\{\mathbf{x}_h\in \mathbf{H}^1_D(\Omega): \mathbf{x}_h|_T\in \mathbf{P}_1(T),\;\; \forall T \in \mathcal{T}_h\}.$ The discrete control space $M_H=\{p_H\in L^2(\Omega): p_H|_T\in P_0(T),\; \forall T\in \mathcal{T}_H\}$. We solve the above discrete optimality \eqref{modified_discModel} system using \textit{primal-dual} active set strategy \cite[Section 2]{Primal_Dual2002}.

 To illustrate the \textit{primal-dual} active set strategy algorithm, let us define some notations. Let the dimension of $\mathbf{V}_h$ and $\mathbf{Q}_h$ be denoted by $2n$ and $2m$, respectively. Also, let the dimension of discrete pressure space $M_H$ to be $\kappa.$ Let $\mathcal{V}_h^D$ and $\mathcal{V}_h^C$ denote the set of vertices on $\overline{\Gamma}_D$ and the set of vertices interior to $\Gamma_C,$ in the fine mesh $\mathcal{T}_h$, respectively. The active and inactive sets corresponding to the bilateral constraints are
\begin{eqnarray}
	\mathcal{A}^k_b &:=& \{ i \in \mathcal{V}_h^C \; | \; \bm{\mu}_i^k + (\mathbf{y}_h^k - \mathbf{y}_b)_i > 0 \}, \nonumber\\
	\mathcal{A}^k_a &:=& \{i \in \mathcal{V}_h^C \; | \; \bm{\mu}_i^k + (\mathbf{y}_h^k - \mathbf{y}_a)_i < 0 \}, \nonumber\\
	\mathcal{I}^k &:=& \{i \in \mathcal{V}_h^C \; | \; \bm{\mu}_i^k + (\mathbf{y}_h^k - \mathbf{y}_b)_i \leq 0 \leq \bm{\mu}_i^k + (\mathbf{y}_h^k - \mathbf{y}_a)_i\}, \nonumber
\end{eqnarray} 
where $\mathbf{y}_h^k$ is the $k^{th}$ iterate of $\mathbf{y}_h,$ and $\bm{\mu} \in \mathbb{R}^{2m}$ is the Lagrange multiplier. We denote by $I$, an identity matrix of size $2m \times 2m$. Moreover, $\mathcal{V}^I_h$ denotes the set of all interior nodes in $\mathcal{T}_h$. Let $\bm{\phi}_i$ are basis functions of $V_h$ and $Q_h.$ Also, let $\chi_i$ are basis function of $M_H.$ We define the following matrices and vectors:
\begin{eqnarray}
	A_{i,j} &:=& \sum_{K \in \mathcal{T}_h} \int_{K} \nabla \bm{\phi}_i:\nabla \bm{\phi}_j,\quad M_{i,j} := \sum_{K \in \mathcal{T}_h} \int_{K} \bm{\phi}_i \cdot \bm{\phi}_j,\quad 	B_{j,k} := \sum_{K \in \mathcal{T}_h} \int_{K} \nabla \cdot \bm{\phi}_j \; \chi_{k}, \nonumber \\
	F^1_j &:=& \sum_{K \in \mathcal{T}_h} \int_{K}  \mathbf{f}\cdot \bm{\phi}_j, \quad F^2_j := \sum_{K \in \mathcal{T}_h} \int_{K}  \mathbf{u}_d \cdot \bm{\phi}_j, \quad L= \rho A+M. \nonumber\\
\end{eqnarray}
\noindent
Now the \textbf{primal-dual active set} algorithm for the Dirichlet boundary control problem \eqref{Modified_prob} reads as: \\
\begin{description}
	\item[Step 1.] Initialize $\mathbf{y}^0_h$ $\bm{\mu}^0$ and set $k = 0$.
	\item[Step 2.] Set the active and inactive sets ($\mathcal{A}^k_a,\mathcal{A}^k_b,\mathcal{I}^k$).
	\item[Step 3.] Solve
	{\small
		\begin{eqnarray}
			A_{2n\times 2n} \left[\mathbf{w}_h^{k+1}\right]_{2n \times 1} + A_{2n \times 2m} \left[\mathbf{y}_h^{k+1}\right]_{2m \times 1}+B_{2n \times \kappa} [p^{k+1}_H]_{\kappa \times 1} &=& [F^1]_{2n \times 1}, \nonumber\\
           B^T_{\kappa\times2n}[\mathbf{w}_h^{k+1}]_{2n\times 1}+B^T_{\kappa\times2m}[\mathbf{y}_h^{k+1}]_{2m\times 1}&=&[\mathbf{0}]_{\kappa\times1}\nonumber\\
			A_{2n \times 2n} \left[\bm{\phi}_h^{k+1}\right]_{2n \times 1}-B_{2n \times \kappa} [r^{k+1}_H]_{\kappa \times 1} - M_{2n \times 2n} \left[\mathbf{w}_h^{k+1}\right]_{2n \times 1} - M_{2n \times 2m} \left[\mathbf{y}_h^{k+1}\right]_{2m \times 1} &=&  -[F^2]_{2n \times 1}, \nonumber\\
           B^T_{\kappa\times 2n} [\bm{\phi}^{k+1}]_{2n\times1} &=& [\mathbf{0}]_{\kappa\times1}\nonumber\\
			L_{2m \times 2m} \left[\mathbf{y}_h^{k+1}\right]_{2m \times 1} - A_{2m \times 2n} \left[\bm{\phi}_h^{k+1} \right]_{2n \times 1}+B_{2m \times \kappa} [r^{k+1}_H]_{\kappa\times1} + M_{2m \times 2n} \left[ \mathbf{w}_h^{k+1}\right]_{2n \times 1} +\nonumber\\  I_{2m \times 2m} \left[\bm{\mu}^{k+1}\right]_{2m \times 1} = [F^2]_{2m\times 1}&+& \rho A_{2m \times 2m} \left[q_d\right]_{2m \times 1}, \nonumber\\
			\left[\mathbf{y}_h^{k+1}\right]_{2m \times 1} &=& \mathbf{y}_a \;\; \mbox{on} \;\; \mathcal{A}^k_a \nonumber\\
			\left[\mathbf{y}_h^{k+1}\right]_{2m \times 1} &=& \mathbf{y}_b \;\; \mbox{on} \;\; \mathcal{A}^k_b, \nonumber\\
			\left[\bm{\mu}^{k+1}\right]_{2m \times 1} &=& \mathbf{0} \;\; \mbox{on} \;\; \mathcal{I}^k \cup \mathcal{V}^I_h. \nonumber
		\end{eqnarray}
	}
	\item[Step 4.] Stop using the criterion $\mathcal{A}^{k+1} = \mathcal{A}^k$ and $\mathcal{I}^{k+1} = \mathcal{I}^k$ or $\norm{\mathbf{y}_h^{k+1}-\mathbf{y}_h^k}_{L^{\infty}(\Omega)} < \epsilon$ for $\epsilon > 0$, or set $k = k + 1$ and return to \textbf{Step 2}.\\
\end{description}

\noindent
In each of our model problems, we compute the error and estimator, which are defined as follows:
\begin{align}\label{error_terms}
  		\textbf{Error}:=\lVert \nabla(\mathbf{y}-\mathbf{y}_h)\rVert_{0,\Omega}+\lVert \nabla(\mathbf{u}-\mathbf{u}_{h})\rVert_{0,\Omega} +\lVert p-p_H\rVert_{0,\Omega}+\lVert \nabla (\bm{\phi}-\bm{\phi}_{h})\rVert_{0,\Omega}+\lVert r-r_H\rVert_{0,\Omega},  
\end{align}
and,
\begin{align}\label{estm_terms}
 \textbf{Estimator} :=\eta_{(\mathbf{u},p)}+\eta_{(\bm{\phi},r)}+\eta_{\mathbf{y}},   
\end{align}
where $\eta_{(\mathbf{u},p)},$ $\eta_{(\bm{\phi},r)}$ and $\eta_{\mathbf{y}}$ are defined in Theorem \ref{aposteriori error estimator}.\\

\noindent
For the adaptive algorithm, we use the following paradigm:
\begin{align*}
   {\it Solve \rightarrow Estimate \rightarrow Mark \rightarrow Refine}.   
 \end{align*}
 
 \noindent
 First, we compute the discrete solutions ($\mathbf{u}_h,p_H,\bm{\phi}_h,r_H,\mathbf{y}_h$) using the above-described primal-dual active set algorithm. Then in the second step using the discrete solution, we compute the error estimator (Estimator = $\eta_{(\mathbf{u},p)}+\eta_{(\bm{\phi},r)}+\eta_{\mathbf{y}}$) over each element. We use the D\"{o}rlfer marking technique {\cite{Dorfler}} with bulk parameter $\theta=0.3$ for the mark step. Then the marked elements are refined using the newest vertex bisection algorithm \cite{adap_bisection} to obtain a new mesh and the algorithm is repeated. The convergence rate for \textbf{Estimator} is defined as follows:
 \begin{align*}
     \text{rate of convergence}(\ell):= \frac{\log(\textbf{Estimator}_{\ell+1}/\textbf{Estimator}_{\ell})}{\log(N_{\ell}/N_{\ell+1})},
 \end{align*}
 for $\ell:= 1,2,3,\cdots$, where $\textbf{Estimator}_{\ell}$ and $N_\ell$ denotes the estimator and number of degrees of freedom at $\ell$th level respectively. Similarly, one can define the rate of convergence for \textbf{Error}.

\begin{figure}
     \centering
     \begin{subfigure}[b]{0.4\textwidth}
         \centering
         \includegraphics[width=\textwidth]{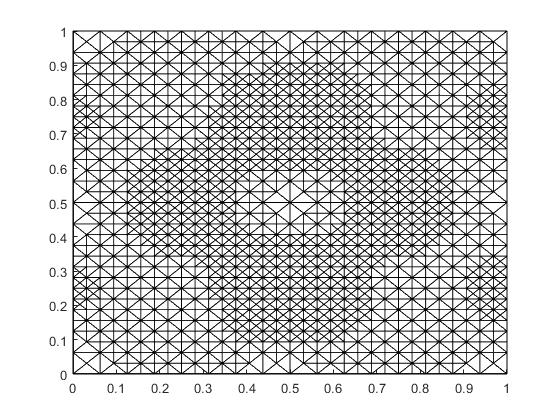}
         \caption{Coarse mesh $\mathcal{T}_H$}
         \label{coarse_mesh}
     \end{subfigure}
     \begin{subfigure}[b]{0.4\textwidth}
         \centering
         \includegraphics[width=\textwidth]{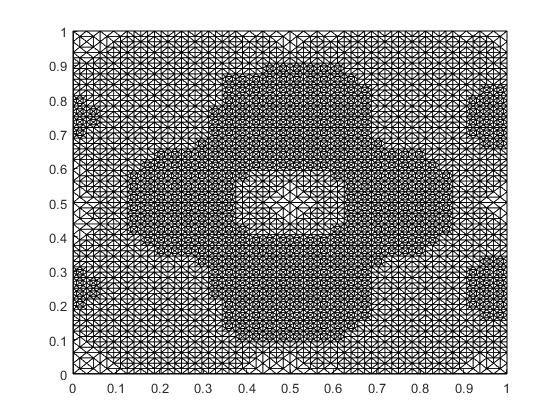}
         \caption{Refined mesh $\mathcal{T}_h$}
         \label{Refine mesh}
     \end{subfigure}
     \begin{subfigure}[b]{0.6\textwidth}
         \centering
         \includegraphics[width=\textwidth]{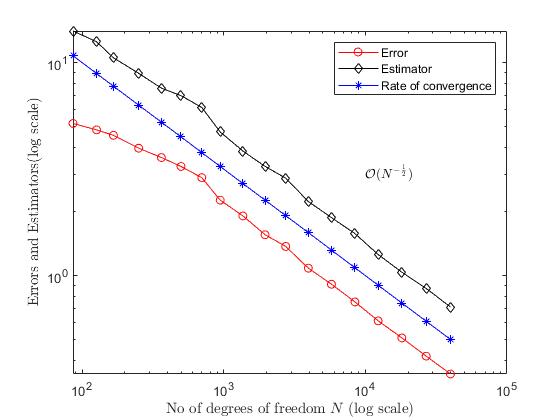}
         \caption{Convergence history (unit square domain)}
         \label{conv_history}
     \end{subfigure}
     \caption{}
     \label{square_dom}
\end{figure}

  \begin{example}\label{ex1}
 	In this example we consider the optimal control problem \eqref{Modified_prob} with the computational domain $\Omega=(0,1)^2$ and $\Gamma_D=(0,1)\times \{0\},$ $\Gamma_C=\partial \Omega \backslash \Gamma_D$. We choose the constants $\rho=10^{-2},$ $\mathbf{y}_a=(-4,-2),$ and $\mathbf{y}_b=(2,2.5).$ The state and adjoint state variables are given by
 	\begin{align*}
 	{\bf u}={\bf y}= \left(
 	\begin{array}{c}
 	-\exp(x)(y\cos(y)+\sin(y))\\
 \exp(x)y\sin(y)
 	\end{array}
 	\right),~~~ p = \sin(2\pi x)\sin(2 \pi y),
 	\end{align*}
 	and
 	\begin{align*}
 	\bm{\phi}= \left(
 	\begin{array}{c}
 	(\sin(\pi x))^2 \sin(\pi y)\cos(\pi y)\\
 	-(\sin(\pi y))^2 \sin(\pi x) \cos(\pi x)
 	\end{array}
 	\right),~~~ r = \sin(2\pi x)\sin(2 \pi y).
 	\end{align*}
 	We choose $\mathbf{u}$ and $\bm{\phi}$ such that $\nabla\cdot{\mathbf{u}}=\nabla\cdot{\bm{\phi}}=0 \quad \text{in}\; \Omega$  and $\bm{\phi}=\mathbf{0}~~ \text{on}~~ \partial\Omega.$ The data of the problem are  chosen such that $\mathbf{f} = -\Delta \mathbf{u}+\nabla{p},~\mathbf{u}_d= \mathbf{u}+\Delta \bm{\phi}+\nabla{r}~ \text{and}~\mathbf{y}_d=\mathbf{y}$.	
 \end{example}
 We have used the above-described primal-dual active set algorithm to solve the optimal control problem. Figure \ref{square_dom}(A) and Figure \ref{square_dom}(B) show the coarse and refine meshes respectively. Figure \ref{square_dom}(C) shows the convergence of error and estimator in terms of the number of degrees of freedom. It is clear from Figure \ref{square_dom}(C), that the \textbf{Estimator} and the \textbf{Error} show the optimal rate of convergence. Here the optimal rate of convergence means the rate of convergence is $0.5$ with respect to the number of degrees of freedom($N$).

\begin{figure}
    \centering
    \includegraphics[width=10cm,height=3.5in]{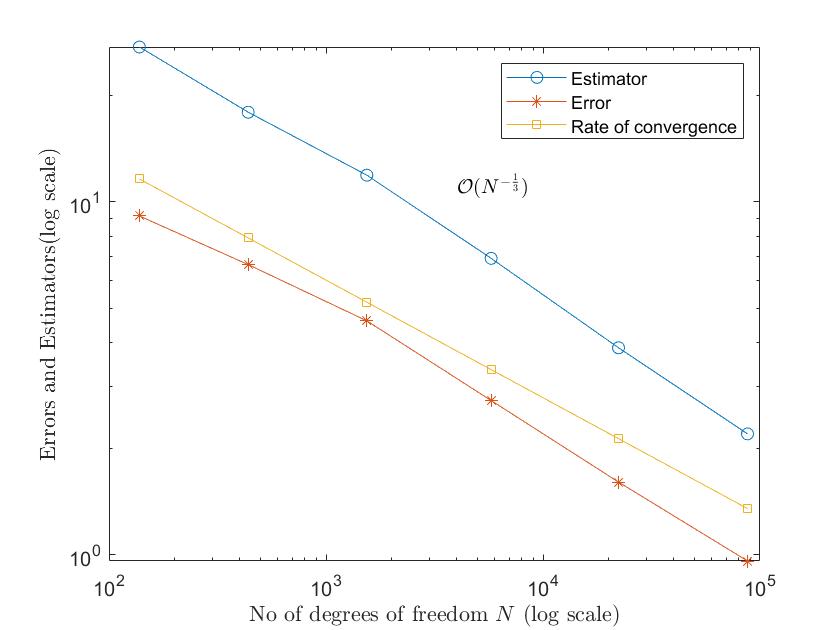}
    \caption{Convergence history on uniform mesh (L-shape domain).}
    \label{convhist_uniform_mesh}
\end{figure}

\begin{figure}
     \centering
     \begin{subfigure}[b]{0.4\textwidth}
         \centering
         \includegraphics[width=\textwidth]{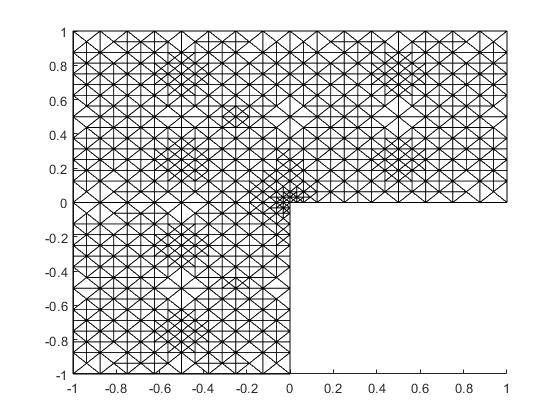}
         \caption{Adaptive coarse mesh $\mathcal{T}_H$}
         \label{Adap:coarse_mesh}
     \end{subfigure}
     \begin{subfigure}[b]{0.4\textwidth}
         \centering
         \includegraphics[width=\textwidth]{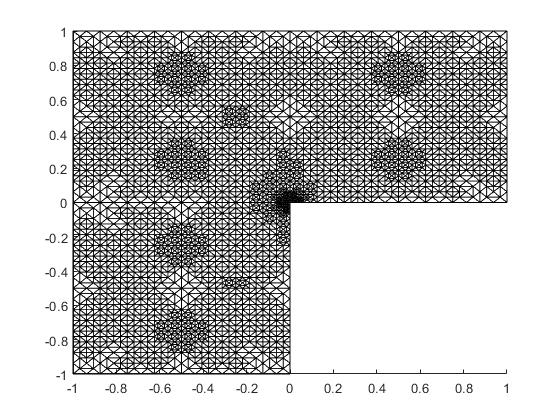}
         \caption{Adaptive refine mesh $\mathcal{T}_h$}
         \label{Adap:Refine mesh}
     \end{subfigure}
     \begin{subfigure}[b]{0.6\textwidth}
         \centering
         \includegraphics[width=\textwidth]{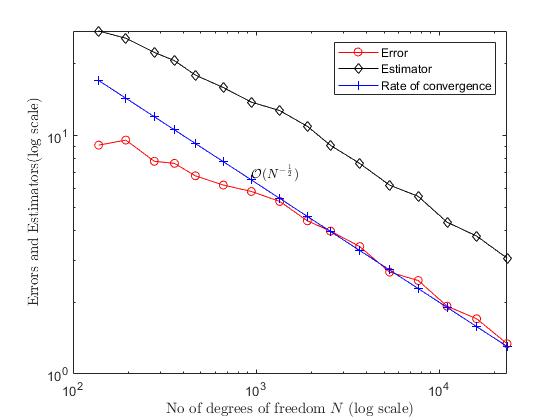}
         \caption{Convergence history(L-shape domain)}
         \label{Adap:conv_history}
     \end{subfigure}
     \caption{}
     \label{LShapeDom}
\end{figure}

 \begin{example}\label{ex2}  In this example we consider the optimal control problem \eqref{Modified_prob} with the  L-shaped domain $\Omega=(-1, 1)^2\setminus ([0,1]\times [-1,0])$, $\Gamma_C=\partial\Omega,$ and the exact solutions
 	\begin{align*}
 		{\bf u}={\bf y} = r^\alpha \left(
 		\begin{array}{c}
 			(1+\alpha)\sin(\theta)\omega(\theta)+\cos(\theta)\omega'(\theta)\\
 			-(1+\alpha)\cos(\theta)\omega(\theta)+\sin(\theta)\omega'(\theta)
 		\end{array}
 		\right),\\
 		p=- r^{\alpha-1}((1+\alpha)^2 \omega'(\theta)+\omega{'''}(\theta))/(1-\alpha),
 	\end{align*}
 	where
 	\begin{align*}
 		\omega(\theta)=&1/(1+\alpha)\sin(\alpha+1)\theta)\cos(\alpha w)-\cos((\alpha+1)\theta)\\ &+1/(1+\alpha)\sin(\alpha-1)\theta)\cos(\alpha\omega)-\cos((\alpha-1)\theta)
 	\end{align*}
 	and $\alpha=856399/1572864$ and $w=3\pi/2$. The adjoint variables ${\bm {\phi}},~r$ are considered the same as in Example \ref{ex1}. The data of the problem is  chosen such that  $\mathbf{f} = -\Delta \mathbf{u}+\nabla{p},~\mathbf{u}_d= \mathbf{u}+\Delta \bm{\phi}+\nabla{r}~ \text{and}~\mathbf{y}_d=\mathbf{y}$.	The constants $\rho = 10^{-2},\; \mathbf{y}_a=(-3,-3),$ and $\mathbf{y}_b=(4,4).$
 \end{example}
 This problem is defined on the L-shaped domain, and the derivative of the solution $({\bf u},p)$ has a singularity at the origin. It is well known that for this problem the uniform refinements will not provide an optimal convergence rate. We have a similar observation from Figure \ref{convhist_uniform_mesh}, for uniform refinements convergence rate with respect to the number of degrees of freedom ($N$) is $0.33$ Hence, one can use the adaptive algorithm to improve the convergence rate. Figure \ref{LShapeDom}(A) and Figure \ref{LShapeDom}(B) show the adaptive coarse and refined meshes respectively. Figure \ref{LShapeDom}(C) shows the adaptive convergence of error and estimator in terms of the number of degrees of freedom($N$). We see that the convergence rate has been improved from $0.33$ (Figure \ref{convhist_uniform_mesh}) to $0.50$ (Figure \ref{LShapeDom}(C)). Thus the optimal convergence is achieved using the adaptive algorithm for the error in energy norm in the state and adjoint state velocity, control approximation, in $L2-$ norm of pressure and adjoint pressure variables. Hence, the optimal convergence for the a posteriori
estimator(\textbf{Estimator}) and the total error(\textbf{Error}), which are defined in \eqref{error_terms} and \eqref{estm_terms}. Here, the optimal convergence means the rate of convergence is $0.5$ with respect to the number of degrees of freedom($N$).

\section{Conclusions}
In this article, we propose, analyze, and test an \textit{a posteriori} error estimator for the Dirichlet boundary control problem governed by Stokes equation. We develop an inf-sup stable finite element discretization scheme by using $\mathbf{P}_1$ elements(in the fine mesh) for the velocity and control variable and $P_0$ elements(in the coarse mesh) for the pressure variable. The optimal control satisfies a bilateral Signorini contact problem, thus the discrete optimality system consists of a discrete variational inequality for the approximate control variable. We derive and analyze the error estimator for the control variable and the estimator is designed for controlling its energy error. The estimator reduces to the standard residual estimator for elliptic problem, if no contact occurs. The contributions by the estimator addressing the nonlinearity are related to the contact stresses and the complementarity condition. 
We prove the reliability and efficiency of the estimator and ensure the equivalence with the error up to oscillation terms. Our numerical experiments confirm the theoretical results.

\bibliography{Control}
\bibliographystyle{abbrvnat}

\end{document}